\documentclass[]{article}
\makeatletter\if@twocolumn\PassOptionsToPackage{switch}{lineno}\else\fi\makeatother

\usepackage{amsfonts,amssymb,amsbsy,latexsym,amsmath,tabulary,graphicx,times,caption,fancyhdr}
\usepackage[utf8]{inputenc}
\usepackage{url,multirow,morefloats,floatflt,cancel,tfrupee}
\makeatletter

\AtBeginDocument{\@ifpackageloaded{textcomp}{}{\usepackage{textcomp}}}
\makeatother
\usepackage{colortbl}
\usepackage{xcolor}
\usepackage{pifont}
\usepackage[nointegrals]{wasysym}
\usepackage[pdftex]{hyperref}

\usepackage[ruled,vlined]{algorithm2e}

\usepackage{float}
\usepackage[utf8]{inputenc}
\usepackage[english]{babel}
\usepackage{graphicx}
\usepackage{caption}
\usepackage{amsthm,amssymb}
\usepackage{amsmath}
\usepackage{thmtools}
\usepackage{mathrsfs}
\usepackage{booktabs}
\usepackage{verbatim}
\usepackage{enumitem}
\usepackage{longtable}
\usepackage{wasysym}
\usepackage{atbegshi}
\usepackage{microtype}
\usepackage{natbib}
\usepackage{setspace}
\usepackage{bbm}

\makeatletter

\def\mcWidth#1{\csname TY@F#1\endcsname+\tabcolsep}

\def\cAlignHack{\rightskip\@flushglue\leftskip\@flushglue\parindent\z@\parfillskip\z@skip}
\def\rAlignHack{\rightskip\z@skip\leftskip\@flushglue \parindent\z@\parfillskip\z@skip}

\@ifundefined{etal}{}{}

\usepackage{ifxetex}
\ifxetex\else\if@twocolumn\@ifpackageloaded{stfloats}{}{\usepackage{dblfloatfix}}\fi\fi

\AtBeginDocument{
\expandafter\ifx\csname eqalign\endcsname\relax
\def\eqalign#1{\null\vcenter{\def\\{\cr}\openup\jot\m@th
  \ialign{\strut$\displaystyle{##}$\hfil&$\displaystyle{{}##}$\hfil
      \crcr#1\crcr}}\,}
\fi
}

\AtBeginDocument{%
  \@ifpackageloaded{endfloat}%
   {\renewcommand\efloat@iwrite[1]{\immediate\expandafter\protected@write\csname efloat@post#1\endcsname{}}}{\newif\ifefloat@tables}%
}%

\def\BreakURLText#1{\@tfor\brk@tempa:=#1\do{\brk@tempa\hskip0pt}}
\let\lt=<
\let\gt=>
\def\processVert{\ifmmode|\else\textbar\fi}

\@ifundefined{subparagraph}{
\def\subparagraph{\@startsection{paragraph}{5}{2\parindent}{0ex plus 0.1ex minus 0.1ex}%
{0ex}{\normalfont\small\itshape}}%
}{}

\newcommand\role[1]{\unskip}
\newcommand\aucollab[1]{\unskip}
  
\@ifundefined{tsGraphicsScaleX}{\gdef\tsGraphicsScaleX{1}}{}
\@ifundefined{tsGraphicsScaleY}{\gdef\tsGraphicsScaleY{.9}}{}
\def\checkGraphicsWidth{\ifdim\Gin@nat@width>\linewidth
	\tsGraphicsScaleX\linewidth\else\Gin@nat@width\fi}

\def\checkGraphicsHeight{\ifdim\Gin@nat@height>.9\textheight
	\tsGraphicsScaleY\textheight\else\Gin@nat@height\fi}

\def\fixFloatSize#1{}
\let\ts@includegraphics\includegraphics

\def\inlinegraphic[#1]#2{{\edef\@tempa{#1}\edef\baseline@shift{\ifx\@tempa\@empty0\else#1\fi}\edef\tempZ{\the\numexpr(\numexpr(\baseline@shift*\f@size/100))}\protect\raisebox{\tempZ pt}{\ts@includegraphics{#2}}}}

\AtBeginDocument{\def\includegraphics{\@ifnextchar[{\ts@includegraphics}{\ts@includegraphics[width=\checkGraphicsWidth,height=\checkGraphicsHeight,keepaspectratio]}}}

\DeclareMathAlphabet{\mathpzc}{OT1}{pzc}{m}{it}

\def\URL#1#2{\@ifundefined{href}{#2}{\href{#1}{#2}}}

\def\UrlOrds{\do\*\do\-\do\~\do\'\do\"\do\-}%
\g@addto@macro{\UrlBreaks}{\UrlOrds}

\edef\fntEncoding{\f@encoding}

\makeatother

\newif\ifmultipleabstract\multipleabstractfalse%
\makeatletter

\def\wileyIndent{1pt}
\usepackage[paperheight=10in,paperwidth=6.5in,margin=2cm,headsep=.5cm,top=2.5cm,headheight=1cm]{geometry}

\renewenvironment{abstract}
{\vspace*{-1pc}\trivlist\item[]\leftskip\wileyIndent\hrulefill\par\vskip4pt\noindent\textbf{\abstractname}\mbox{\null}\\}{\par\noindent\hrulefill\endtrivlist}

\usepackage[]{footmisc}

\def\author#1{\gdef\@author{\hskip-\dimexpr(\tabcolsep)\hskip\wileyIndent\parbox{\dimexpr\textwidth-\wileyIndent}{\centering\bfseries#1}}}

\def\title#1{\linespread{1}\gdef\@title{\centering\bfseries\ifx\@articleType\@empty\else\@articleType\\\fi#1}}

\let\@articleType\@empty \def\articletype#1{\gdef\@articleType{{\normalfont\itshape#1}}}

\AtBeginDocument{\fancypagestyle{headings}{\fancyhf{}\fancyhead[C]{\RunningHead}\fancyhead[R]{\thepage}}\pagestyle{headings}}

 \def\audegree#1{}

\declaretheoremstyle[%
  spaceabove=12pt,%
  spacebelow=12pt,%
  headfont=\normalfont\itshape,%
  postheadspace=1em,%
]{definition}

\declaretheoremstyle[%
  spaceabove=12pt,%
  spacebelow=12pt,%
  headfont=\normalfont\itshape,%
  postheadspace=1em,%
]{thm} 
\declaretheorem[name={\bf Theorem},style=thm,
]{thm}

\declaretheorem[name={\bf Lemma},style=thm,
]{lemma}

\declaretheorem[name={\bf Corollary},style=thm,
]{cor}

\declaretheoremstyle[%
  spaceabove=-6pt,%
  spacebelow=6pt,%
  headfont=\normalfont\itshape,%
  postheadspace=1em,%
  qed=\qedsymbol%
]{pf} 
\declaretheorem[name={Proof},style=pf,unnumbered,
]{pf}

\declaretheoremstyle[%
  spaceabove=24pt,%
  spacebelow=6pt,%
  headfont=\normalfont\itshape,%
  postheadspace=1em,%
]{remark}

\newcommand{\Mod}[1]{\ (\mathrm{mod}\ #1)}

\usepackage{authblk}

\date{}

\title{The Number of Ribbon Tilings for Strips}
\author[1]{Yinsong Chen}
\author[2]{Vladislav Kargin}
\affil[1]{Department of Mathematics and Statistics, Binghamton University \protect\\ Binghamton, New York, U.S.A.
\protect\\ ychen276@Binghamton.edu}
\affil[2]{Department of Mathematics and Statistics, Binghamton University
\protect\\  Binghamton, New York, U.S.A 
\protect\\ vkargin@Binghamton.edu}
\def\RunningHead{}

\begin{document}

\maketitle 

\begin{abstract}
First, we consider order-$n$ ribbon tilings of an $M$-by-$N$ rectangle $R_{M,N}$ where $M$ and $N$ are much larger than $n$. We prove the existence of the growth rate $\gamma_n$ of the number of tilings and show that $\gamma_n \leq (n-1) \ln 2$. Then, we study a rectangle $R_{M,N}$ with fixed width $M=n$, called a strip. We derive lower and upper bounds on the growth rate $\mu_n$ for strips as $	\ln n - 1 + o(1) \leq \mu_n \leq \ln n $. Besides, we construct a recursive system which enables us to enumerate the order-$n$ ribbon tilings of a strip for all $n \leq 8$ and calculate the corresponding generating functions.

\def\keywordstitle{Keywords}

\end{abstract}

\section{Introduction}
\begin{figure}[H] 
	\centering
	\includegraphics[scale=0.35]{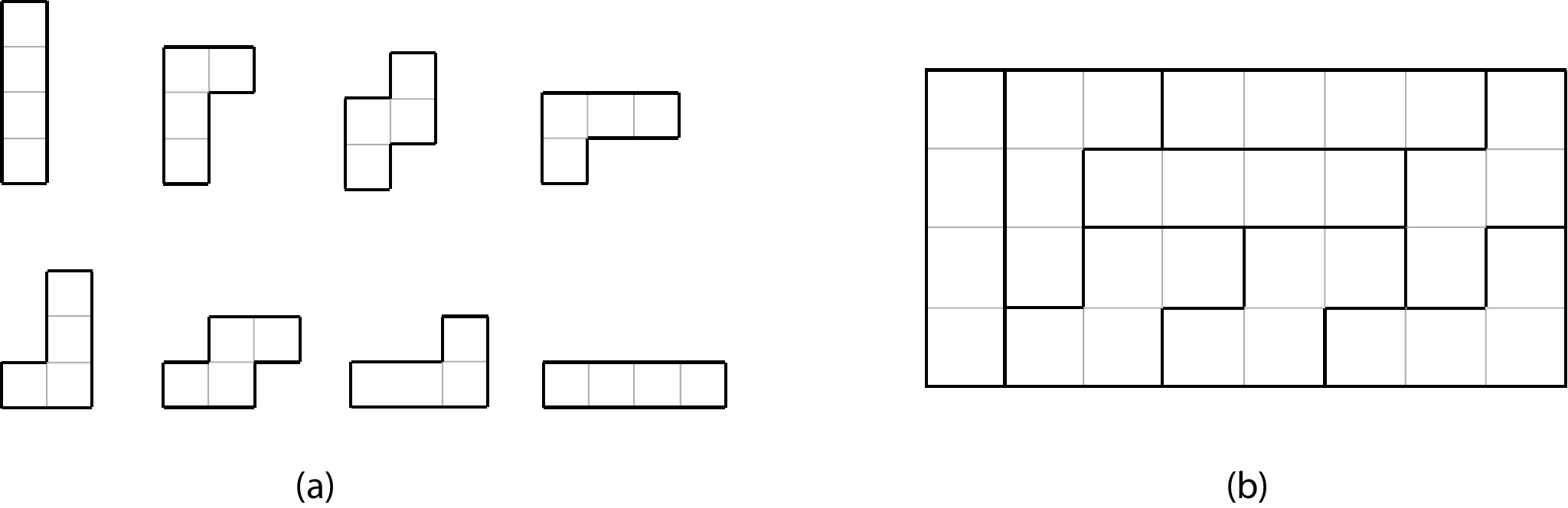} 
	\caption{For the order $n=4$, (a) eight types of ribbons (the picture is from \citep{sheffield2002ribbon}); and (b) a tiling of a strip with dimension $4 \times 8$.}  
	\label{FigTilingOrder4}
\end{figure}

Let an integer $n \geq 2$ be fixed. We say a square $[x, x+1] \times [y, y+1]$ in the two-dimensional integer lattice has \emph{level} $l = x+ y$ and \emph{color} $c \equiv x+y  \pmod{n}$ where $0 \leq c < n$. An order-$n$ \emph{ribbon tile} is a set of squares connected along an edge and containing exactly one square of each color. In other words, a ribbon tile is a sequence of $n$ adjacent squares, each of which is located above or to the right of its predecessor. An order-$n$ \emph{ribbon tiling} is a covering of a region with non-overlapping order-$n$ ribbon tiles. (See Figure \ref{FigTilingOrder4} for an example of ribbon tiles and a ribbon tiling when $n=4$.)

Two natural questions about ribbon tilings are whether there exists a ribbon tiling of a given region, and if it does, then how many different ribbon tilings there are. The existence question for every simply connected region and arbitrary order $n$ was settled in \citep{sheffield2002ribbon}, who provided an algorithm, linear in the area of a region, that checks if the region has a ribbon tiling. The enumeration question for domino tilings, which are ribbon tilings for $n=2$, goes back to 1960's. The papers \citep{kasteleyn1961statistics} and  \citep{temperley1961dimer} provided a formula for rectangular regions. \citep{klarner1980domino} and \citep{stanley1985dimer} studied domino tilings of rectangular regions with fixed width from a different perspective. In particular, they studied the generating function for the number of domino tilings. Later, \citep{kenyon2016asymptotics} considered domino tilings on a torus. A technique to analyze domino tilings of more general regions was invented by \citep{conway1990tiling} and \citep{thurston1990conway}. A related enumeration problem of tilings with T-tetrominoes was studied by \citep{korn2004tilings}. Ribbon tilings for $n > 2$ were first studied in  \citep{pak2000ribbon}. Pak introduced the term \emph{tile counting group} and made a conjecture about ribbon tiling, which was proved in \citep{moore2002ribbon}. These techniques have been extended in \citep{sheffield2002ribbon} who discovered a one-to-one correspondence between ribbon tilings and acyclic orientations of a certain partially oriented graph. This discovery led him to the development of the algorithm verifying the existence of ribbon tilings. \citep{alexandersson2018enumeration} calculated the number of ribbon tilings for $n>2$ of a special rectangle with dimension $n \times 2n$. In this paper, we extend these results by studying the number of ribbon tilings of rectangular regions for $n > 2$. Specially, we study the $n \times N$ strips with arbitrary large $N$. For $N = 2n$, our numerical results are in agreement with the results in \citep{alexandersson2018enumeration}.

\subsection{Main results}
Let $R_{M,N}$ be an $M$-by-$N$ rectangle. It is known that $R_{M,N}$ has an order-$n$ ribbon tiling if and only if $n | M$ or $n | N$. In this paper, we always suppose $n | M$. Let $f_{M,N}$ be the number of ribbon tilings of $R_{M,N}$. Define the \emph{growth rate} of $f_{M,N}$ to be 
\begin{equation}		\label{EqGrowthRateDef}
	\gamma_n = \lim_{M,N \to \infty} \frac{ \ln (f_{M,N} )}{ A }
\end{equation}
where $A = MN / n$ is the number of tiles of $R_{M,N}$, provided the limit in Equation \ref{EqGrowthRateDef} exists. Our first result is demonstrating the existence of the growth rate $\gamma_n$.

\begin{thm}  \label{ThmGrowthRateExist}
	For any order $n \geq 2$, the growth rate $\gamma_n$ exists and $\gamma_n \leq (n-1) \ln 2$. 
\end{thm}

The proof of Theorem \ref{ThmGrowthRateExist} will be given in Section 2. 

Now, we consider a special rectangle $R_N$ with fixed width $n$ and length $N$. We call $R_N$ a \emph{strip} with length $N$. It is evident that $R_N$ can be tiled by $N$ rectangles with dimension $n \times 1$. Hence, the region $R_N$ always has at least one order-$n$ ribbon tiling. We are interested in the \emph{growth rate} of the number $f_N$ of ribbon tilings of $R_N$, that is 
$$
	\mu_n = \lim_{N \to \infty} \frac{\ln f_N}{N}.
$$
The existence of $\mu_n$ can be proved by an argument similar to that in the proof of Theorem \ref{ThmGrowthRateExist}. We give upper and lower bounds on $\mu_n$ in Theorem \ref{ThmGrowthRate}.

\begin{thm}  \label{ThmGrowthRate}
	The growth rate $\mu_n$ satisfies the following inequalities:
	$$ 
	\frac{1}{n} \ln (n!) \leq \mu_n \leq \ln (n).
	$$
\end{thm}

\begin{cor} \label{CorGrowthRate} For growing $n$, 
	$$ 
	\ln (n) - 1 + o(1) \leq \mu_n \leq \ln (n).
	$$
\end{cor}

The proof of Theorem \ref{ThmGrowthRate} and Corollary \ref{CorGrowthRate} will be given in Section 4.

For the number of tilings $f_N$, we define \emph{the leftmost tiling process} in order to produce a recursive system. In the leftmost tiling process, there are $(n-1)!$ specific types of regions with different left boundaries that re-appear repeatedly in the process. We call these $(n-1)!$ types the \emph{fundamental regions}. One of these fundamental regions, with a vertical line as the left boundary, has the same shape as the original $n$-by-$N$ strip. Let  $\mathbf{f}(N)$ be a vector of which each component is the number of ribbon tilings of the corresponding fundamental region with size $N-i$ where $0 \leq i < n$. Theorem \ref{ThmDiffEq} shows that the vector $\mathbf{f}(N)$ satisfies a recursive system. 

\begin{thm} \label{ThmDiffEq}
Consider order-$n$ ribbon tilings of an $n$-by-$N$ strip. The number of ribbon tilings of fundamental regions satisfies a recursive system 
\begin{equation*}		\label{EqRecurrence}
	\mathbf{f}(N) = A_n \mathbf{f}(N-1)	.
\end{equation*}
The transfer matrix $A_n$ has the following properties:
\begin{itemize}
	\item[(1)] The elements of $A_n$ are either $0$ or $\pm 1$.
	\item[(2)] $A_n$ has the form
$$
A_n = 
  \left[
  \begin{matrix}
   A_n' & I_2  \\
   I_1 & 0  
  \end{matrix} 
  \right] .
$$
where $A_n '$ has dimension $(n-1)! \times (n-1)(n-1)!$, $I_1$ is an identity matrix with dimension $(n-1)(n-1)! \times (n-1)(n-1)!$ and $I_2$ is an identity matrix with dimension $(n-1)! \times (n-1)!$.
\end{itemize}
\end{thm}

\begin{figure}[H]
	\centering
	\includegraphics[scale=0.55]{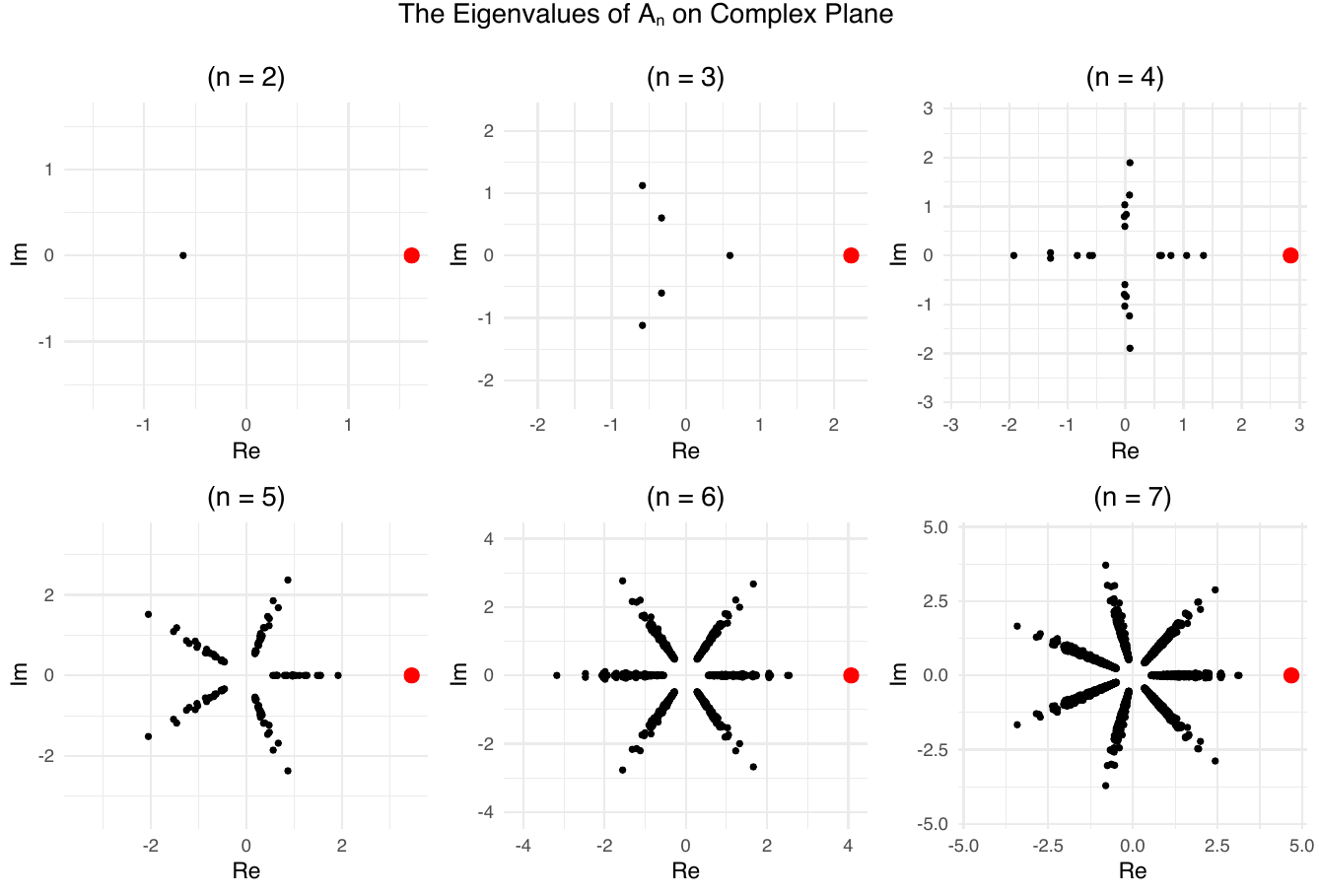}
	\caption{The plot shows the eigenvalues of $A_n$ in the complex plane for $n = 2, 3, 4, 5, 6,7$. The eigenvalue with the largest absolute value (the growth rate) is depicted with a red dot.}
	\label{FigEigenValuesDist}
\end{figure}

The construction of the recursive system and the proof of Theorem \ref{ThmDiffEq} are in Section 5. There is some interesting structure in the eigenvalues of $A_n$ shown in Figure \ref{FigEigenValuesDist}.

\begin{figure}[H]
	\centering
	\includegraphics[scale=0.55]{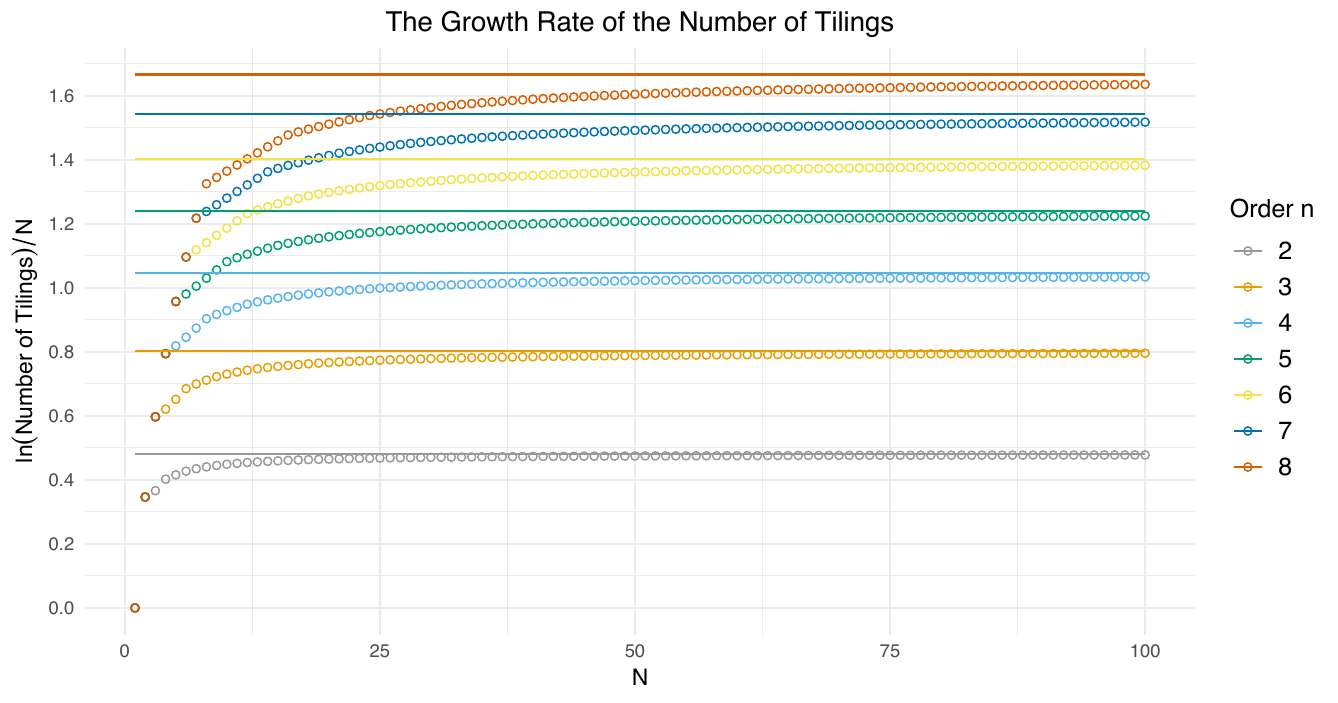}
	\caption{For order $n=2, 3,\cdots,8$, the plot shows the asymptotic behavior of the number of ribbon tilings of a strip. The horizontal lines are the growth rates $\mu_n = \ln (\lambda_n)$. }
	\label{FigGrowthRate}
\end{figure}

\begin{figure}[H]
	\centering
	\includegraphics[scale=0.5]{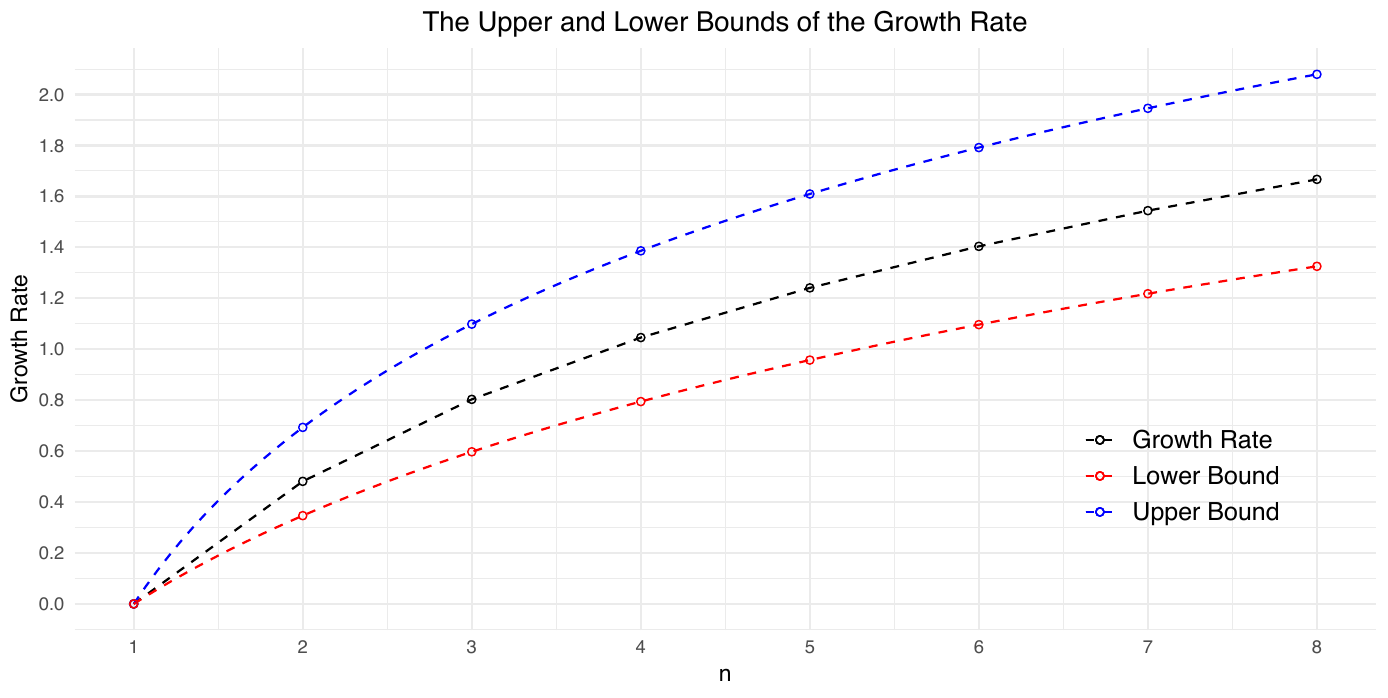}
	\caption{For order $n=2, 3,\cdots,8$, the plot shows the upper and lower bounds on the growth rates given by $\frac{1}{n} \ln (n!) \leq \mu_n \leq \ln (n)$ in Theorem \ref{ThmGrowthRate}.}
	\label{FigGrowthRateBound}
\end{figure}

Numeric results for the number of ribbon tilings $f_N$ complement our result in Theorem \ref{ThmGrowthRate}. Table \ref{TableNumTiling} lists the first ten terms of $f_N$ for the cases $n = 2, 3,\cdots, 8$. Let $\lambda_n$ be the eigenvalue with the largest absolute value of the transfer matrix $A_n$ in Theorem \ref{ThmDiffEq}. Then, $\mu_n = \ln(\lambda_n)$. Numeric values of $\lambda_n$'s are listed in Table \ref{TableGrowthRate}; they are rounded to six digits after decimal. Figure \ref{FigGrowthRate} shows the asymptotic behavior of $f_N$. Figure \ref{FigGrowthRateBound} compares the growth rates with their bounds in Theorem \ref{ThmGrowthRate}.

\begin{table}[H] 
\centering
\begin{tabular}{rrrrrrrrrrrr}
  \hline
 N & 1 & 2 & 3 & 4 & 5 & 6 & 7 & 8 & 9 & 10 & $\cdots$ \\ 
  \hline
  n=2 & 1 & 2 & 3 & 5 & 8 & 13 & 21 & 34 & 55 & 89 & $\cdots$ \\ 
  n=3 & 1 & 2 & 6 & 12 & 26 & 61 & 134 & 297 & 669 & 1490 & $\cdots$ \\ 
  n=4 & 1 & 2 & 6 & 24 & 60 & 160 & 455 & 1379 & 3849 & 10811 & $\cdots$ \\ 
  n=5 & 1 & 2 & 6 & 24 & 120 & 360 & 1140 & 3810 & 13434 & 49946 & $\cdots$ \\ 
  n=6 & 1 & 2 & 6 & 24 & 120 & 720 & 2520 & 9240 & 35490 & 142758 & $\cdots$ \\ 
  n=7 & 1 & 2 & 6 & 24 & 120 & 720 & 5040 & 20160 & 84000 & 364560 & $\cdots$ \\ 
  n=8 & 1 & 2 & 6 & 24 & 120 & 720 & 5040 & 40320 & 181440 & 846720 & $\cdots$ \\ 
   \hline
\end{tabular}
\caption{The numbers of ribbon tilings of an $n$-by-$N$ strip by order-$n$ ribbons.} 
\label{TableNumTiling}
\end{table}

\begin{table}[H] 
\centering
\resizebox{\columnwidth}{!}{
\begin{tabular}{cccccccc}
\hline
                         & $n=2$     & $n=3$      & $n=4$     & $n=5$     & $n=6$ & $n=7$ & $n=8$  \\
\hline
$\lambda_n$	& $1.618034$     & $2.232476$ & $2.845807$ & $3.458663$ & $4.071277$ & $4.683752$ & $5.296141$ \\
\hline
\end{tabular}
}
\caption{The eigenvalues with the largest absolute value of transfer matrices.} 
\label{TableGrowthRate}
\end{table}

The rest of the paper is organized as follows. In Section 2, we use the superadditivity argument to show the existence of the growth rate $\gamma_n$. In Section 3, we introduce the leftmost tiling process in order to enumerate all ribbon tilings of a strip. In Section 4, we prove Theorem \ref{ThmGrowthRate}. In Section 5, we construct the recursive system to calculate the number of tilings of a strip.

\section{Superadditivity}
Each pair of tilings of non-overlapping regions $R_1$ and $R_2$ corresponds to a tiling of the region $R_1 \cup R_2$. From this fact, it follows that the logarithm of the number of tilings is superadditive. We formalize this argument in Lemma \ref{LemSupadd}. Lemma \ref{LemGrowRateExist} is a multidimensional version of Fekete's lemma. In all statments, it is assumed that $n | M$.

\begin{lemma} 	\label{LemUpBoundSimple}
Let $f_{M,N}$ be the number of tilings in an $M$-by-$N$ rectangle. Then, 
	$$
		f_{M,N} \leq 2^{\frac{MN(n-1)}{n}}	.
	$$
\end{lemma}

\begin{pf}
	There are $M N / n$ tiles in an $M$-by-$N$ rectangle. By using the result of \citep{sheffield2002ribbon}, they can be canonically ordered. Let the tiles be labeled as $1, 2, \cdots, M N / n$, which is invariant for all tilings. In different tilings, each labeled tile may have different types. By the definition of a ribbon tile, there are $2^{n-1}$ different types for a tile. Hence, a tiling corresponds to a sequence of types, which consist of $MN/n$ elements and this map is injective. It is not necessary for each sequence of types to be a valid tiling. Then, the number of tilings can be upper bounded as
	$$
		f_{M,N} \leq (2^{n-1})^{MN/n} = 2^{\frac{MN(n-1)}{n}}	.
	$$
\end{pf}

\begin{lemma}	\label{LemSupadd}
	Suppose $M_2 = p M_1$ and $N_2 = q N_1$. Then, $\ln f_{M_2, N_2} \geq p q \ln f_{M_1, N_1}$.
\end{lemma}

\begin{pf}
	A $M_2$-by-$N_2$ rectangle can be separated  into $pq$ number of $M_1$-by-$N_1$ rectangles. The union of the tilings of small rectangles is a valid tiling of a big rectangle. This map is injective. 
\end{pf}

\begin{lemma}	\label{LemGrowRateExist}
The following equality holds
	$$ 
		\lim_{M, N \to \infty} \frac{\ln ( f_{M,N} ) }{MN} = \sup_{M,N \geq 1} \frac{\ln ( f_{M,N} ) }{MN} < \infty.
	$$
\end{lemma}

\begin{pf}
	Let $S = \sup_{M,N \geq 1} \frac{\ln (f_{M,N})}{MN}$. Then $S$ is finite by Lemma \ref{LemUpBoundSimple}. For any $\epsilon > 0$, the definition of $S$ gives $M_0$ and $N_0$ such that $ \frac{\ln (f_{M_0,N_0})}{M_0 N_0} >  (S-\epsilon)$. Consider an $M$-by-$N$ rectangle $R$ and pick integers $p, q$ such that $p M_0 \leq M \leq (p+1) M_0$ and $q N_0 \leq N \leq (q+1) N_0$. Since every ribbon tiling of the left lower $p M_0$-by-$q N_0$ sub-rectangle of $R$ can be extended to a ribbon tiling of $R$, we get $f_{M,N} \geq f_{p M_0, q N_0}$. Then Lemma \ref{LemSupadd} gives 
	
$$
	\frac{\ln f_{M,N}}{MN}  \geq \frac{\ln f_{p M_0, q N_0}}{MN}  \geq \frac{p q \ln f_{M_0,N_0}}{(p+1)M_0 (q+1)N_0}
	\geq \frac{p}{p+1} \frac{q}{q+1} (S - \epsilon)  .
$$	
	
Since $\epsilon$ is arbitrary, we obtain $ \lim_{M,N \to \infty} \frac{\ln ( f_{M,N} ) }{MN} = S$.
\end{pf}

Now, we are ready to prove Theorem \ref{ThmGrowthRateExist}.
\begin{pf}
Lemma \ref{LemGrowRateExist} establishes that the growth rate $\gamma_n$ exists. Lemma \ref{LemUpBoundSimple} implies the upper bound on $\gamma_n$.
\end{pf}

\section{The Leftmost Tiling Process}
First, let us explain how ribbon tilings are related to acyclic orientations on partially oriented graphs. The paper \citep{sheffield2002ribbon} introduced a ``\emph{left of}'' relation for both tiles and squares, denoted as $\prec$. Let $s_{x,y}$ be a square $[x, x+1] \times [y, y+1]$. We say $s_{x,y} \prec s_{x',y'}$ if one of the following two conditions holds:
\begin{itemize}
	\item[(1)] $x + y = x' + y'$ and $x < x'$;
	\item[(2)] $|(x+y)-(x'+y')| = 1$, $x \leq x'$ and $y \geq y'$.
\end{itemize}

Let $t$ be a tile and $s$ be a square. We write $s \prec t$ if $s \prec s'$ for some square $s' \in t$, and $t \prec s$ if $s' \prec s$ for some square $s' \in t$. If $t_1$ and $t_2$ are two tiles in a tiling, we write $t_1 \prec t_2$ if there exist a square $s_1 \in t_1$ and a square $s_2 \in t_2$ with $s_1 \prec s_2$. It is not possible that both $t_1 \prec t_2$ and $t_2 \prec t_1$ unless $t_1 = t_2$. The relation $\prec$ is not transitive. However, if $t_1 \prec t_2 \prec \ldots \prec t_k$, then either $t_1 \prec t_k$  or $t_1$ and $t_k$ are incomparable. It cannot happen that $t_k \prec t_1$.

We say that a tile has \emph{level} $l$ if $l$ is the lowest level of the squares in this tile. \citep{sheffield2002ribbon} showed that every tiling of a region $R$ has the same number of tiles in a given level. Hence tiles can be enumerated independently of a tiling as ``a tile number $i$ in level $l$''. Note that in different tilings a tile with the same label can have different locations and different shapes. Given this enumeration, we identify tiles with vertices in a graph $G_R$. The set of vertices of $G_R$ consists of all tiles and boundary squares (squares outside $R$, but having an edge in $R$), and two vertices are connected by an edge if they are comparable with respect to the ``left of'' relation. Further, Sheffield showed that there is a one-to-one correspondence between ribbon tilings and acyclic orientations of the graph $G_R$ with a fixed partial orientation.

\begin{figure}[H]
			\centering
			\includegraphics[scale=0.55]{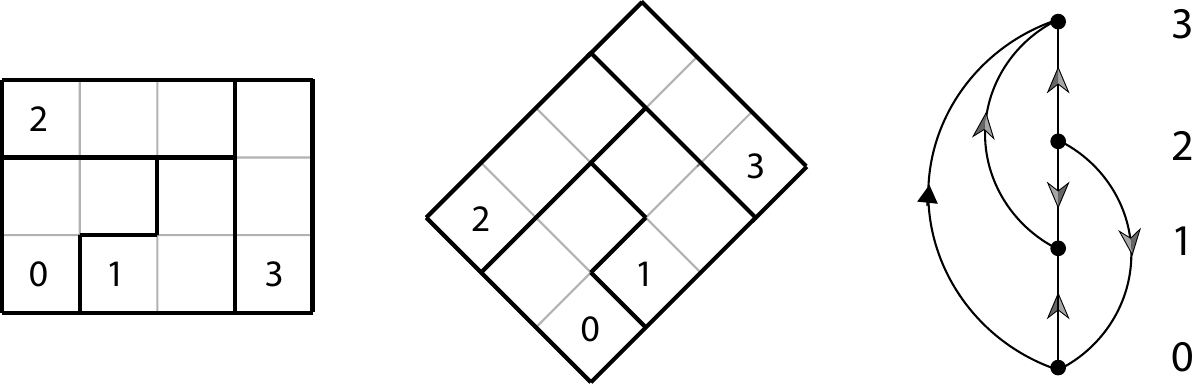}
			\caption{$n=3$. 	Observe the ``left of'' relation by rotating the picture $45$ degrees. In the corresponding acyclic orientation, note that $0 \to 3$ is forced. The tile with label $2$ is the leftmost tile in this tiling.} 
			\label{FigCorresTilingAO}
	\end{figure}

For an $n$-by-$N$ strip region $R$, the graph $G_R$ can be described as follows. Note that the ``left of'' relations between boundary squares and tiles are the same in all tilings and a boundary square can never be between two tiles, so we can restrict $G_R$ to the vertices that represent tiles. Sheffield’s results imply that there is exactly one tile in each level for a strip $R$. Therefore, we label the vertices of $G_R$ with the level of their corresponding tiles as $0, 1, \cdots, N-1$. There is an edge between two vertices $i$ and $j$ if and only if $|i-j| \leq n$. The edge is oriented from $i$ to $j$, denoted as $i \to j$, if $i \prec j$. Note that the orientation $i \to j$ is forced if $j-i=n$ since in all tilings the ``left of'' relation $i \prec j$ is fixed for $j - i = n$. This fixed orientation on the edge $(i, j)$ is a part of the fixed partial orientation on $G_R$, and these are the only forced orientations on $G_R$. The other orientations depend on
the particular tiling of $R$. See Figure \ref{FigCorresTilingAO} for an example of the ``left of'' relation and the acyclic orientation correspondence.

Let $0, 1, \cdots, N-1$ be the labels of tiles that will be used to cover an $n$-by-$N$ strip $R$. It was shown by \citep{sheffield2002ribbon} that building a tiling is equivalent to determining the ``left of" relations among these labeled tiles. We will build these relations step by step using the concept of the ``leftmost'' tile. Let $\alpha$ be a tiling of $R$ and $G_R(\alpha)$ be the acyclic orientation of $G_R$ corresponding to $\alpha$. Every tile that corresponds to a source of $G_R(\alpha)$ is called a \emph{source tile}. The \emph{leftmost tile} is a source tile with the smallest label. 


Now, we introduce \emph{the leftmost tiling process} to determine relations among tiles step by step. At each step, we choose an appropriate label and declare that the tile with this label will be the \emph{leftmost} tile in a tiling of the untiled region. We put this tile in place and continue the procedure. We keep choosing a label from the set of remaining labels and declaring the corresponding tile the leftmost among the remaining tiles until all the relations among all tiles are determined. 

We say a sequence $T$ of tile labels is a \emph{tiling sequence} if it is a valid sequence for the leftmost tiling process. It is evident that there is a bijection between the set of ribbon tilings and the set of tiling sequences. We can enumerate ribbon tilings of $R$ with the help of the leftmost tiling process. In Figure \ref{FigCorresTilingAO}, the tiling corresponds to the tiling sequence $T = [2, 0, 1, 3]$. 

\begin{figure}[H]
	\centering
	\includegraphics[scale=0.35]{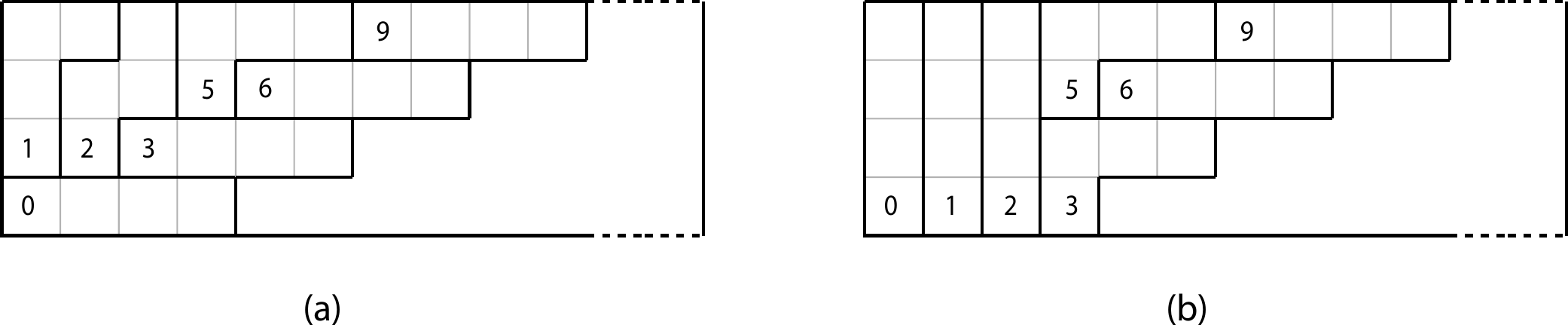}
	\caption{$n=4$. The tiling sequence $T = [ 1, 2, 5, 9, 6, 3, 0 ]$ is shown in (a), while the tiling sequence $\zeta_0(T) = [ 0, 1, 2, 5, 9, 6, 3 ]$ is shown in (b). The residual region shown in (a) and (b) are the same.}
	\label{FigReturnOperator}
\end{figure}

Let us define an operation on tiling sequences, which will be useful later. For any sequence 
	$$
		T = [ t_1, \cdots, t_k, 0, s_1, \cdots, s_l]
	$$ 
containing a unique $0$, define the \emph{return operator} $\zeta_0$ by the formula 
	$$
		\zeta_0 (T) = [ 0, t_1, \cdots, t_k,  s_1, \cdots, s_l].
	$$
That is, the operator $\zeta_0$ moves the $0$ to the front of the sequence. See Figure \ref{FigReturnOperator} for an example of the return operator. This operator will be used in the proof of both Theorem \ref{ThmGrowthRate} and Theorem \ref{ThmDiffEq}.

\begin{lemma} \label{LemValidSeq0}
	Let $T$ be a tiling sequence of a strip $R$. Then $\zeta_0 (T)$ is a valid tiling sequence of $R$.
\end{lemma}

\begin{pf}
	Let $G_R$ be the partially ordered graph associated to the rectangle $R$.
	Let $T = [ t_1, \cdots, t_k, 0, t_{k+2}, \cdots, t_N ]$ correspond to a tiling of $R$ with associated orientation $\alpha$ on $G_R$. If $G_R(\alpha)$ has edges between some of $t_1, \cdots, t_k$ and $0$, they are directed towards $0$. Otherwise, $0$ would be chosen before one of these $t_1, \cdots, t_k$ in the tiling sequence. By reversing directions of these edges, we define another orientation $\beta$ of $G_R$. This reversal does not affect the forced orientations on $G_R$, because no forced edge is directed towards $0$. Moreover, $G_R(\beta)$ is acyclic since $0$ is a source in $G_R(\beta)$, so there are no cycles that go through $0$, and cycles that do not include $0$ do not exist, too, because $G_R(\alpha)$ is acyclic. So $\beta$ represents a ribbon tiling of $R$. To see that $\zeta_0(T)$ is the tiling sequence of that tiling, note first that $0$ is a source in $G_R(\beta)$. Moreover, since $t_s$, $1 \leq s \leq k$, was a source in $G_R(\alpha) \setminus \{ t_1, \cdots, t_{s-1} \}$, it is a source in $G_R(\beta) \setminus \{ 0, t_1, \cdots, t_{s-1} \}$. So the new sequence can start with $[0; t_1, \cdots, t_k]$. The final part $ [ t_{k+2}, \cdots, t_N ] $ of the tiling sequence of $G_R(\alpha)$ agrees with that of $G_R(\beta)$, because $G_R(\alpha) \setminus \{ t_1, \cdots, t_k, 0 \} = G_R(\beta) \setminus \{ t_1, \cdots, t_k, 0 \}$.
\end{pf}

\section{Upper and Lower Bounds on the Number of Tilings}
In this section, we prove Theorem \ref{ThmGrowthRate} and Corollary \ref{CorGrowthRate}. For a fixed integer $n \geq 2$, let $f_N$ be the number of order-$n$ ribbon tilings of a strip with length $N$. First, the function $\ln(f_N)$ is superadditive, that is, $\ln f_N \geq \ln f_{N_1} + \ln f_{N_2}$ where $N = N_1 + N_2$. This is because every pair of tilings of two non-overlapping regions $R_1$ and $R_2$ forms a tiling of $R_1 \cup R_2$. By Fekete's Lemma, we have the existence of the following limit
$$ 
\lim_{N \to \infty} \frac{\ln f_N}{N} = \mu_n 
$$  
where $\mu_n = \sup_N \{ \frac{\ln f_N}{N} \} $ could be a constant or infinity. 
We call $\mu_n$ the \emph{growth rate} of the number of ribbon tilings of a strip.  In order to prove the upper bound on $\mu_n$ stated in Theorem \ref{ThmGrowthRate}, we first prove the following lemma.

\begin{lemma} \label{LemMaxNum}
Let $f_N$ be the number of ribbon tilings of a strip $R$ with length $N$. Then, 
	$$
		f_N \leq n f_{N-1} .
	$$ 
\end{lemma}

\begin{pf}
Let $\mathbb{T}$ be the set of all tiling sequences of $R$ and $\mathbb{T}_0$ be the set of all tiling sequences of $R$ starting with $0$. Note that $f_N = | \mathbb{T} |$ and $f_{N-1} = | \mathbb{T}_0 | $ by the one-to-one correspondence between tilings and tiling sequences. We need to show that $ | \mathbb{T} | \leq n | \mathbb{T}_0 | $. 

By Lemma \ref{LemValidSeq0}, for every $T \in \mathbb{T}$, $\zeta_0(T)$ is a valid tiling sequence of $R$. It follows that the image set $\zeta_0(\mathbb{T})$ is a subset of $\mathbb{T}_0$. Define $C_i = \{ T \in \mathbb{T}: \zeta_0(T) = T_i \}$ where $\{ T_i \}$ is an enumeration of $\zeta_0(\mathbb{T})$. It is clear that $\{ C_i \}$ are mutually exclusive and $\cup C_i = \mathbb{T}$. 

We claim that the size of each $C_i$ is less than or equal to $n$. Indeed, in order to recover a tiling sequence $T \in \zeta_0^{-1} (T_i) = C_i$ from $T_i \in \zeta_0(\mathbb{T}_0)$, it is enough to consider the possible positions at which the tile $0$ can be embedded into $T_i$. By the definition of the leftmost tiling process, either tile $0$ is the first element, or it follows a comparable tile $t$. Otherwise, if $t > 0$ and $0$ were not comparable, then $0$ would have been declared the leftmost before $t$. It follows that the possible tiles that can be followed by $0$ are $1,2,\cdots,n-1$. Then, there are $n$ possible positions for $0$ to be embedded into $T_i$ to recover a tiling sequence $T \in C_i$. Hence, $| C_i | \leq n$ for all $i$. Thus, we have $| \mathbb{T} |   \leq n | \mathbb{T}_0 |$.
\end{pf}

Now, we are ready to prove Theorem \ref{ThmGrowthRate} and Corollary \ref{CorGrowthRate}.
\begin{pf}
For the upper bound, by using the inequality in Lemma \ref{LemMaxNum}, we obtain 
$$
	f_N \leq n f_{N-1} \leq n^2 f_{N-2} \leq \cdots \leq n^N .
$$
It follows that $\mu_n = \lim_{N \to \infty} \frac{\ln f_N}{N} \leq \ln n$. 

For the lower bound, consider a strip $R_m$ that has length $m$ for $m \leq n$. The corresponding graph of $R_m$ is the complete graph with $m$ vertices and no forced edges. Then, tilings of $R_m$ bijectively correspond to permutations of $\{ 0, 1, \cdots, m-1 \}$. It follows that the number of tilings $f_m$ of $R_m$ equals $m!$. Next, we can divide an $n$-by-$N$  strip $R$ into $\lfloor \frac{N}{n} \rfloor$ squares with dimensions $n \times n$ and a small remainder rectangle with dimension $n \times m$ where $m = N \Mod{n}$. By superadditivity, we have the inequality
	$$
		\ln f_N \geq \lfloor \frac{N}{n} \rfloor \ln (n!) + \ln m!	
	$$
By the definition of $\mu_n$, we obtain the lower bound in Theorem \ref{TableGrowthRate}.

By Stirling's approximation $\ln(n!) = n \ln (n) - n + O( \ln(n) )$, we obtain Corollary \ref{CorGrowthRate} from Theorem \ref{ThmGrowthRate}.
\end{pf}

\section{Enumeration}
This section explains the enumeration process in detail and proves Theorem \ref{ThmDiffEq}. 

At each step of the leftmost tiling process, it is important to know which tiles are valid candidates for the leftmost tile. Lemma \ref{LemCandMin} allows one to find all valid choices for the leftmost tile at each step of the process. 

\begin{figure}[H]
	\centering
	\includegraphics[scale=0.4]{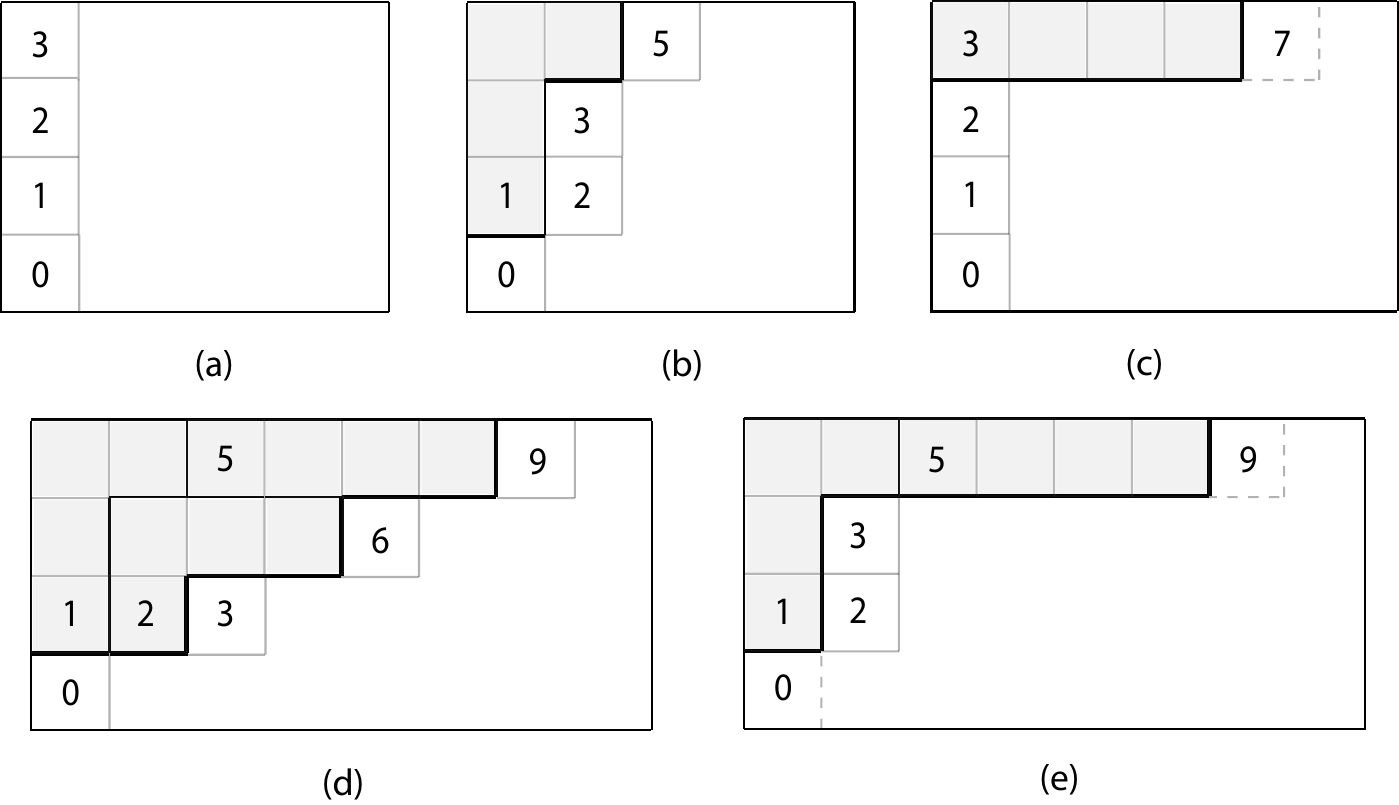}
	\caption{$n=4$. The grey area is the tiled region and the bold lines show the boundary of the current tiling. The squares with grey line boundary indicate the possible candidates for the leftmost location. The squares with dashed line boundary indicate some of those that are ruled out by Lemma \ref{LemCandMin}.} 
	\label{FigCandMin}
\end{figure}

The leftmost tiling process builds a sequence of labels $S=[ w_1, w_2, \cdots, w_N ]$ as an output. At each step of the process, let $X$ be the sequence of labels that have been declared at the previous steps and let $Y = S \backslash X$ be the set of remaining labels. Let $w_X \in X$ be the final element of $X$ and $\max(X)$ be the largest element of $X$. Note that $\max(X)$ can be different from $w_X$. For example, in Figure \ref{FigCandMin} (d), $X= [ 1, 5, 2 ]$, so $\max(X)=5$ and $w_X = 2$. Obviously, $w_X \leq \max(X)$.

Let $\bar{C} = \{ w \in Y :  w_X - n < w \leq \max(X) + n \text{ and } w - n \notin Y\}$ in the case when $X$ is not empty, and let $\bar{C} = \{ w \in Y: 0 \leq w \leq n-1 \}$ when $X$ is empty (that is, at the first step). Sort $\bar{C}$ as an increasing sequence  $[ w_1,w_2, \cdots, w_k ]$. Let $i$ be the smallest index such that $w_{i+1} - w_i \geq n$ for $1 \leq i \leq k-1$ where $k$ is the largest index of the sequence $\bar{C}$. Let $C=\{ w_1, \cdots, w_i \}$ if $i$ exists, and $C=\bar{C}$ otherwise. 

\begin{lemma} \label{LemCandMin}
	At each step of the leftmost tiling process, $w_c \in Y$ can be the leftmost tile of a tiling of the remaining region if and only if $w_c \in C$.
\end{lemma}

Initially, the tile covering the north-western corner of the strip is always the leftmost tile. It follows that the candidate set is $\{ 0, 1, \cdots, n-1\}$ at the first step. The complete proof of Lemma \ref{LemCandMin} is based on the study of the corresponding acyclic orientation of $G_R$ and is relegated to Appendix. In order to give an intuitive idea, we illustrate the lemma for $n=4$ in an example shown in Figure \ref{FigCandMin}. 

Diagram (a) shows that the candidate labels are $0,1,2,3$ at the first step of the leftmost tiling process. Lemma \ref{LemCandMin} says that $4$ is not a valid candidate because $4 - n = 0 \in Y$. (Intuitively, at the first step there are no tilings such that 4 is a source tile with the smallest label.) 

Diagram (b) shows that the candidate labels are $0,2,3,5$ at the second step if $1$ was declared at the first step. Note that tile $5$ is a valid candidate since $5-n=1 \in X$ and so this choice is not ruled out by Lemma \ref{LemCandMin}. (Intuitively, if one declares 5 a source tile with the smallest label, one can build a tiling such that this declaration is true.) 

In Diagram (c), we see that $7$ is not a valid candidate since $7 - 2 = 5 > n$, so that $7$ is in $\bar{C}$ but not in $C$. Again, intuitively this is because $7$ can never be a leftmost tile in the tiling of the remaining region. Diagram (d) shows a case that $C = \bar{C} = \{ 0, 3, 6, 9 \}$. 

In Diagram (e), $0$ is not a valid candidate since $5 - 0 = 5 > n$, so the first inequality in the definition of $\bar{C}$ is not satisfied. Also note that $9$ is in $\bar{C}$ but not in $C$. (Intuitively, 0 is not valid candidate in diagram (e) because declaring it the leftmost tile would contradict the choice of 5 as the source tile with the smallest label in the previous step.) 

Lemma \ref{LemCandMin} provides us an efficient way to find out all the leftmost tiles at each step. We can enumerate all tiling sequences using Lemma \ref{LemCandMin}.


\subsection{The Deducting Moment and Fundamental Regions}
In the leftmost tiling process of a strip $R$, we can stop at any specific step and obtain a sequence $J$. We call $J$ \emph{initial segment} of a tiling sequence. The region $R$ is separated into two parts: a partial region covered by the tiling corresponding to $J$ and a remainder region that is not yet tiled. Define a \emph{residual region} $R_J$ to be the untiled region which is obtained from $R$ by removing the tiling corresponding to $J$. The region $R_J$ has $N-|J|$ tiles. See Figure \ref{FigReturnOperator} for an example. 

More generally, we use the notation $R_J(N-l)$ for $0 \leq l \leq N$ to denote a residual region of the strip with length $N - l + |J| $. That is, $R_J(N-l)$ contains $N-l$ tiles and its left boundary is a horizontal translate of that of $R_J$. 
Clearly, $R_J(N-|J|) = R_J$ since $R_J$ contains $N-|J|$ tiles. For any two residual regions $R_J(N-l)$ and $R_{J'}(N-l')$, we say that they are \emph{similar} if their left boundaries are translates of each other. Obviously, $R_J(N-l)$ is similar to $R_J$ for any valid $N-l \in \mathbb{N}$. The next lemma will show that a residual region $R_J$ is invariant under any valid permutation of $J$. 

\begin{lemma} 	\label{LemPermInvariant}
Let $J$ be an initial segment. If a permutation $\pi(J)$ of $J$ is also a valid initial segment, then $R_J = R_{\pi(J)}$.
\end{lemma}

\begin{pf}
Let $Y$ be the label set of tiles in the region $R_J$. For the graph $G_R$, consider the cut set of $(J, Y)$. The two sequences $J$ and $\pi(J)$ give the same orientation for this cut set, in which each edge is directed from $J$ (or $\pi(J)$) to $Y$. So we can use the same orientations on the subgraph induced by $Y$ and be sure that they will lead to acyclic orientations whether we add them to orientations defined by $J$ or to orientations defined by $\pi(J)$. In particular, we can use orientations on $Y$ defined by a sequence $[ w_{|J| + 1}, \cdots, w_N ]$. Then we can build the tiling starting the construction from the end of this sequence and proceeding by going backwards to the beginning and adding new tiles. This method will determine the identical tilings on $R_J$ and $R_{\pi(J)}$. In particular, this shows that these regions coincide.
\end{pf}

In the leftmost tiling process, define the \emph{deducting moment} to be the step at which the tile with label $0$ (which is unique) is declared to be the leftmost. At the deducting moment, we call a column of the strip a \emph{full column} if it is fully covered by the current tiling. We will show later in Lemma \ref{LemMaxStep} that the deducting moment will be achieved in $O(n^2/2)$ steps. 

\begin{figure}[H]
	\centering
	\includegraphics[scale=0.35]{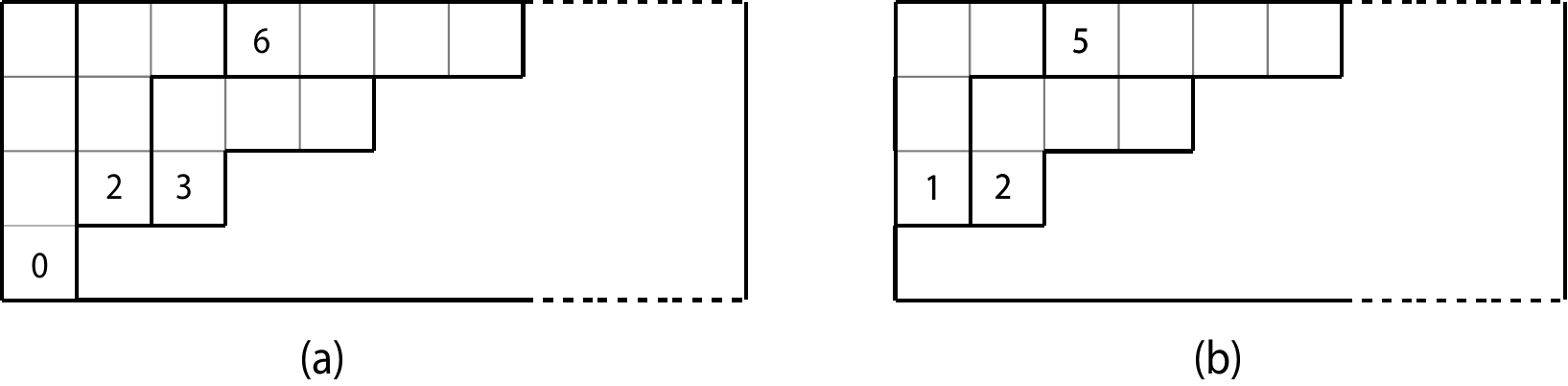}
	\caption{$n=4$. $J = [0, 2, 6, 3]$ and $\delta(J) = [1,5, 2]$. The tiling in (b) is obtained from that in (a) by removing the first column and relabeling the tiles. } 
	\label{FigDeductionOperator}
\end{figure}

For any sequence 
	$$
		J = [0, j_1, \cdots, j_k]
	$$ 
starting with $0$, define the \emph{deduction operator} $\delta$ by the formula 
	$$
		\delta(J) = [ j_1 - 1, \cdots, j_k -1 ] 	.
	$$ 
That is, the operator $\delta$ removes the zero from the sequence and decreases each element by $1$. For an initial segment $J$ starting with $0$, $\delta$ can be regarded as an operator that removes the fist column of a strip and relabel the tiles. That is, $\delta(J)$ is a tiling sequence of a strip with length $N-1$ and its corresponding ribbon tiling is the same as that of $J$ with the first column removed. See Figure \ref{FigDeductionOperator} for an example.

A residual region obtained at the deducting moment can be deducted to a similar residual region by removing initial full columns. Let $J$ be an initial segment achieving the deducting moment with $d$ full columns. Then $0$ is the last element of $J$, and $0, 1, \cdots, d-1 \in J$. We use the composite operator $\delta \circ \zeta_0$ repeatedly $d$ times and denote it as $(\delta \circ \zeta_0)^d$. Let $J' = (\delta \circ \zeta_0)^d (J)$. The deducted region $R_{J'}$ is similar to $R_J$ and their size difference is $d$.

Now, we are ready to define \emph{fundamental regions}. For an initial tiling sequence $J$ achieving the deducting moment with $d$ full columns, we call the sequence 
\begin{equation}		\label{EqFundSeqDef}
	F = \phi (J) = \pi( (\delta \circ \zeta_0 )^d (J) )
\end{equation}
a \emph{fundamental sequence}, where $\pi$ by definition permutes every sequence as an increasing sequence. Let $\mathbb{F}$ be the set of all fundamental sequences, and $R_F(N-l)$ be a \emph{fundamental region} corresponding to $F \in \mathbb{F}$ with size $N-l$. Note that the empty set is a valid fundamental sequence. For notational convenience, we use $[0]$ to represent the empty fundamental sequence since their corresponding fundamental regions are similar. By convention, we will keep $\mathbb{F}$ ordered first by the increasing length of the fundamental sequences and then by the lexicographic order if the sequences have the same length.  

Algorithm \ref{AlgGenFundSeq} is provided to obtain all fundamental sequences by running the leftmost tiling process. We list the results for $n=3$ and $n=4$ as examples. For $n=3$, the fundamental sequences are $[0]$ and $[1]$. For $n=4$, the fundamental sequences are $[0]$, $[1]$, $[2]$, $[1,2]$, $[1,3]$ and $[1,2,5]$. 

Every fundamental sequence $F$ is obtained from an initial segment $J$ ending with $0$ and having $d \geq 1$ full columns. Namely, $F = \pi( (\delta \circ \zeta_0 )^d (J) ) =  \delta^d ( \pi(J) )$. Hence, it will be useful to discuss the properties of $\pi(J)$. We will show that $\pi(J)$ is a valid initial segment in Lemma \ref{LemValidIncreaseSeq}. This implies that a fundamental sequence $F =  \delta^d ( \pi(J) )$ is a valid initial segment. In order to prepare for the proof of Lemma \ref{LemValidIncreaseSeq}, we introduce an auxiliary sequence $A_J$.

Let $A_J = \{w \in J: w+n \notin J \}$. By the definition of the leftmost tiling process, it is not possible that $w+n \in J$ if $w \notin J$. It follows that $A_J$ is a subset of $J$ obtained by keeping the largest element of each modulo-$n$ equivalence class. We permute $A_J = [a_1, \cdots, a_k]$ as an increasing sequence and call it the  \emph{essential sequence} of $J$. Since $J$ ends with $0$, it follows that $m n \notin J$ for any positive integer $m$. Thus, we have $0 \in A_J$ and $a_1 = 0$. (For example, in Figure \ref{FigTranFundamentaSeq1} (a), the initial sequence $[1,2,5,3,0]$ has the essential sequence $[0,2,3,5]$.) 

We say that an increasing sequence $W = [w_1, \cdots, w_k]$ has \emph{connectivity} if $w_{i+1} - w_i < n$ for $i=1,2, \cdots, k-1$. Note that this is stronger than the theoretical connectivity of subgraph $[w_1, \cdots, w_k]$ in $G_R$: we rule out $w_{i+1} - w_i =n$ although $w_i$ and $w_{i+1}$ are connected by an edge. The reason for our definition is that: only if $W$ has connectivity, then it is possible to have a directed path $w_k \to w_{k-1} \to \cdots \to w_1$ without generating any directed cycle in the graph $G_R$. 

Since $J$ is an initial sequence ending with $0$, for every $w \in J \setminus \{0\}$, we must have a directed path from $w$ to $0$ corresponding to the tiling decided by $J$. From the existence of these paths, it follows that $\pi(J)$ has connectivity. Otherwise, let $\pi(J) = [j_1=0, \cdots, j_s, j_{s+1}, \cdots, j_k]$ such that $j_{s+1} - j_s \geq n$, then it is not possible for the directed edge $j_{s+1} \to j_s$. Together with the fact that there is no edge between $j_{s+1}$ and $j_1, \cdots, j_{s-1}$, it follows that it is not possible to have a directed path from $j_{s+1}$ to $j_1$ in the graph $G_J$. This contradicts the assumption that $J$ is an initial segment ending with $0$. 

Lemma \ref{LemGenSeqConnectivity} shows that the essential sequence $A_J$ also has connectivity.

\begin{lemma}	\label{LemGenSeqConnectivity}
	Let $J$ be an initial segment ending with $0$ and let $A_J$ be the essential sequence of $J$. Then, $A_J$ has connectivity.
\end{lemma}

\begin{proof}
	In order to seek contradiction, suppose $a_{m-1}, a_m \in A_J$ such that $a_m - a_{m-1} > n$. (It is not possible $a_m - a_{m-1} = n$ by the definition of $A_J$.) Let $D_0 = \{ w \in J: a_m - n \leq w < a_m \}$ and $D_1 = \{ w+n: w \in D_0 \}$. By our assumption, $D_0 \cap A_J = \emptyset$. It follows that if $w \in D_0$, then $w+n \in J$, otherwise $w \in J$ would be the largest in its equivalence class and would belong to $A_J$. Thus, we have $D_1 \subseteq J$. 

	From the definition of $D_0$ and $D_1$, all elements of $D_1$ are comparable, thus $G_{D_1}$ is a complete graph. Let $s$ be the unique sink of $G_{D_1}$ corresponding to the orientation induced by the initial segment $J$. We have already noted that for every $w \in J$ there is a directed path from $w$ to $0$. In particular, it follows that there is a directed path $p_s$ in $G_J$ from $s$ to $0$. 
	
	Since $D_0$ contains all elements of $J$ between $a_m-n$ and $a_m$, it follows that the directed path $p_s$ must contain a directed edge $s \to w_0$ for some $w_0 \in D_0$. However, $s \to w_0 \to w_0+n \to s$ is a directed cycle where $w_0+n \in D_1$, contradiction.
\end{proof}

Now, we are ready to prove Lemma \ref{LemValidIncreaseSeq}.

\begin{lemma}	\label{LemValidIncreaseSeq}
	Let $J$ be an initial segment ending with $0$ and let $\pi(J)$ be the increasing permutation of $J$. Then, $\pi(J)$ is a valid initial segment.
\end{lemma}

\begin{proof}
	For the increasing sequence $\pi(J)$, $0$ is declared to be the leftmost tile at the first step. From Lemma \ref{LemCandMin}, it is clear that $0$ is a valid candidate at the first step. 
	
	Let $u$ be the $l$-th element and $v$ be the $(l+1)$-th element of the sequence $\pi(J)$, and $\pi(J)_u$ be the sub-sequence of $\pi(J)$ containing all elements before and including $u$. We will prove the lemma by induction. Suppose $\pi(J)_u$ is a valid initial segment. We need to show that $v$ is a valid candidate for the leftmost tile at the $(l+1)$-th step.

	At the $(l+1)$-th step, we have $X = \pi(J)_u$ and $w_X = \max(X) = u$. Recall that $\bar{C} = \{ w \in Y :  w_X - n < w \leq \max(X) + n \text{ and } w - n \notin Y\}$. Then, in our situation, we have $\bar{C} = \{ w \in Y  :  u - n < w \leq u + n \text{ and } w - n \notin Y \}$. 
	
	Let $A_J$ be the essential sequence of $J$, and let $A_u$ be the sub-sequence of $A_J$ such that $A_u = \{ w \in A_J: u-2n < w \leq u \}$, then for every $w \in A_u$ we have $u-n < w+n \leq u+n$. By the definition of essential sequence, it follows that for every $w \in A_u$ we have $w+n \notin J$. Let $A'_u = \{w+n: w\in A_u \}$. By the definition of $\bar{C}$, we have $A'_u \subseteq \bar{C}$. 
	
	Let $w_0$ be the smallest element of $\bar{C}$. It is clear that $w_0 > u-n$ by the definition of $\bar{C}$. By Lemma \ref{LemGenSeqConnectivity}, $A_J$ has connectivity and $0 \in A_J$, it follows that $\min(A_u) - (u-2n) < n$ by the definition of $A_u$, and thus $\min(A'_u) - (u-n) = \min(A_u) + n - (u-n) < n$. Therefore, we have $\min(A'_u) - w_0 < n$. Note that the connectivity of $A'_u$ is inherited from $A_u$. It follows that $A'_u \subseteq C$ by the definition of $C$.
	
	By the definition of $\bar{C}$, it is clear that $u+n \in \bar{C}$. From the connectivity of $A_J$, it follows that $0 \leq u - \max(A_u) < n$ by the definition of $A_u$, and thus $(u+n) - \max(A'_u) < n$. Since $\max(A'_u) \in C$ and $u+n \in \bar{C}$, we have $u+n \in C$ by the definition of $C$. Note that $u+n$ is the largest element of $\bar{C}$. Therefore, we have $C = \bar{C}$.
	
	From the connectivity of $\pi(J)$, it follows that $u < v < u+n$. By the definition of $\bar{C}$, it follows that $v \in \bar{C}$, and thus $v \in C$ since $C = \bar{C}$. Therefore, $v \in \pi(J)$ is a valid candidate at the $(l+1)$-th step. By induction, $\pi(J)$ is a valid initial segment.
\end{proof}

Lemma \ref{LemValidIncreaseSeq} shows that every fundamental sequence $F$ is a valid initial segment. Let $C^0_F$ be the candidate set at the first step of the leftmost tiling process in a fundamental region $R_F$. We obtain $C^0_F$ by setting $w_X = 0$ and $X=F$ in Lemma \ref{LemCandMin}. Note that the order of $F$ does not make any difference for $C^0_F$ except through setting $w_X = 0$, since the fundamental region $R_F$ is invariant under any valid permutation of $F$ by Lemma \ref{LemPermInvariant}.

Note, however, that in the leftmost tiling process of the strip $R$, if $J$ is an initial segment, the order of $J$ (and in particular, the last element of $J$) plays an important role in the determination of the candidate set $C_J$ of the leftmost tiles by Lemma \ref{LemCandMin}. In particular, if $F$ is a fundamental sequence and $R_F$ is similar to  $R_J$ then it might happen that $C_J \neq  C^0_F$, which will lead to an issue for our  enumeration method. This issue will be discussed and rectified in the next section.

\subsection{The Recursive System}
We think about the tiling process as the sequence of transitions between fundamental regions (up to similarity). By running the leftmost tiling process in a fundamental region $R_F$, we will obtain for the first time another fundamental region in one of the following two situations. (One of the two situations must happen.) Let $J$ be a sequence obtained by running the leftmost tiling process in $R_F$. Then, 
\begin{itemize}
	\item Case (1): $J$ achieves the deducting moment with $d$ full columns.
	\item Case (2): $J$ does not achieve the deducting moment, but $\pi ([F, J]) $ is a fundamental sequence.
\end{itemize}
We consider each of these cases as a transfer between two fundamental sequences. Define $F \xrightarrow{J} F'$ to be a \emph{tiling transition} from $F$ to $F'$
such that 
	$$
		F' = \left\{
		\begin{aligned}
			& \phi ( [F, J]) & \text{ if Case (1)}; 	\\
			& \pi ([F, J]) & \text{ if Case (2)}.
		\end{aligned}
		\right.
	$$
where $\phi$ is as defined in Equation \ref{EqFundSeqDef}. 
For a tiling transition $F \xrightarrow{J} F'$, we call $J$ a \emph{transition sequence}. Note that $R_{F'}(N-l-|J|)$ is the first fundamental region hit by the leftmost tiling process in the region $R_{F}(N-l)$. (We assume here that $N-l$ is not too small so that the region $R_{F}(N-l)$ can be tiled by $J$.) 

\begin{figure}[H]
	\centering
	\includegraphics[scale=0.35]{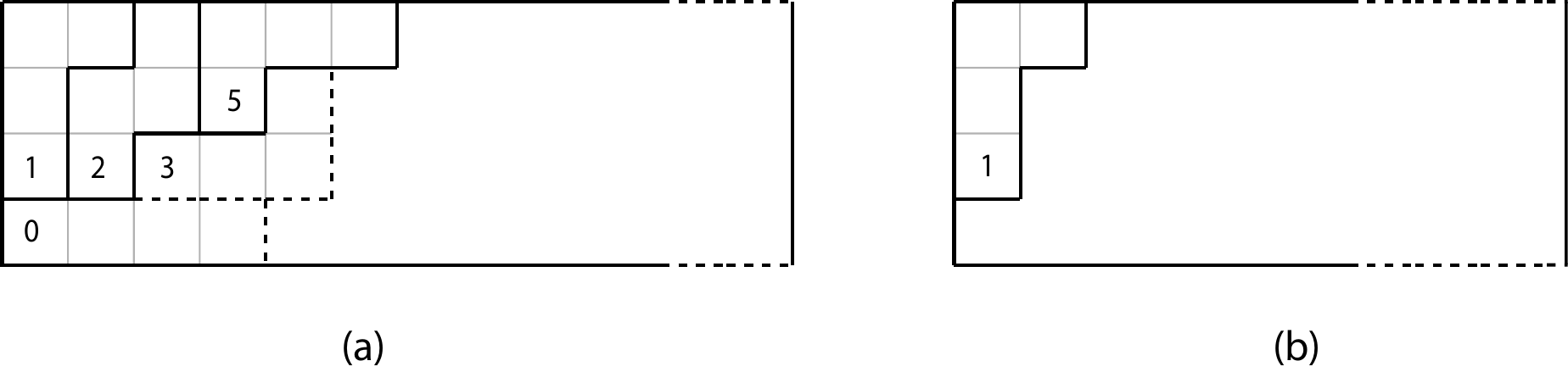}
	\caption{$n=4$. Given the region $R_F$ with $F = [ 1,2,5 ]$ shown in (a) with black bold boundary and the transition sequence $J = [ 3,0 ]$, we obtain $F' = \phi ( [F, J] ) = \phi ( [ 1,2,5,3,0 ] ) = [ 1 ]$ and the corresponding region is shown in (b).} 
	\label{FigTranFundamentaSeq1}
\end{figure}

\begin{figure}[H]
	\centering
	\includegraphics[scale=0.35]{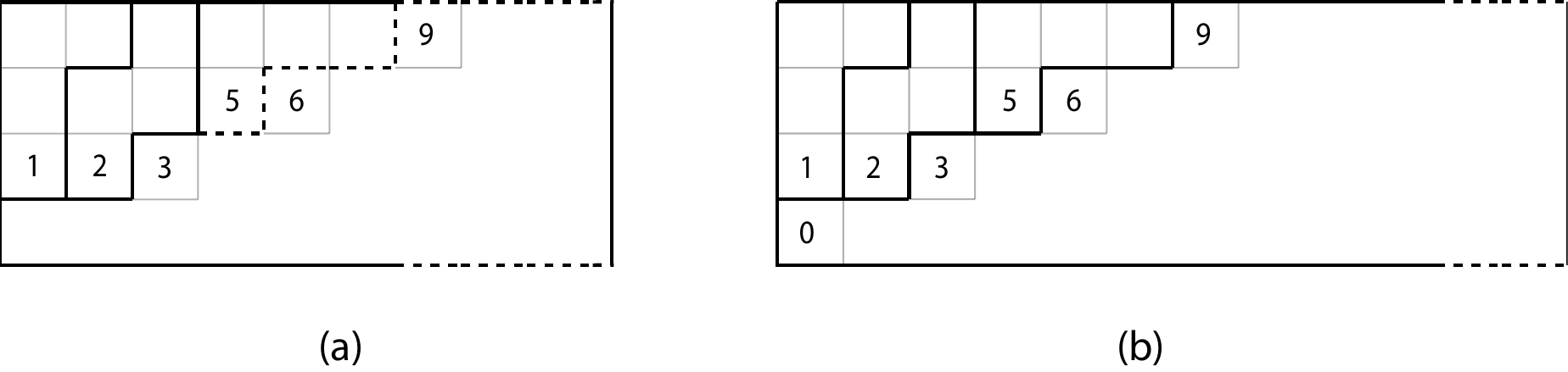}
	\caption{$n=4$. Given the region $R_F$ with $F = [1, 2]$ shown in (a) with black bold boundary and the transition sequence $J = [ 5 ]$, we obtain $F' = \pi ( [F, J] ) = [ 1, 2, 5 ]$ and the corresponding region is shown in (b).} 
	\label{FigTranFundamentaSeq2}
\end{figure}

Figure \ref{FigTranFundamentaSeq1} shows an example of a Case (1) tiling transition. Suppose we start from the region $R_F$ with $F = [ 1,2,5 ]$ shown in (a) with black bold boundary. The transition sequence is $J = [ 3,0 ]$. Then, we obtain $F' = \phi ( [F, J] ) = \phi ( [ 1,2,5,3,0 ] ) = [ 1 ]$ and the corresponding region is shown in (b). In Figure \ref{FigTranFundamentaSeq2}, we show an example of a Case (2) tiling transition from $F=[1,2]$ to $F'=[1,2,5]$ with transition sequence $J=[5]$.

Let $f_F(N-j)$ be the number of tilings of the fundamental region $R_F(N-j)$. In order to calculate $f_F(N-|F|)$, we run the leftmost tiling process in $R_F$ and consider all possible tiling transitions $F \xrightarrow{J} F'$. By the definition of fundamental sequence, $R_{F'}$ and $R_{[F,J]}$ are similar regions, hence $R_{F'}(N-|F|-|J|)$ and $R_{[F,J]}(N-|F| - |J|)$ are identical. Therefore, it is reasonable to collect $f_{F'}(N-|F|-|J|)$ as a positive term in the equation calculating $f_F(N-|F|)$.

	By collecting $f_{F'}(N-|F|-|J|)$ in the equation, we are trying to use $f_{F'}(N-|F|-|J|)$ to represent the number of tiling sequences of $R_F$ starting with $J$. This might lead to an error since $f_{F'}(N-|F|-|J|)$ is obtained by ignoring the order of the sequence $[F,J]$, as mentioned in the last paragraph in Section 5.1. 
	
	Taking the summation over all possible tiling transitions $\tau: F \xrightarrow{J} F'$ in the fundamental region $R_F$, we have 

\begin{equation}		\label{EqNumTilingTran}
	f_F(N-|F|) = \sum_{\tau} f_{F'}(N-|F|-|J|) + e
\end{equation}
where $e$ is an error term.

Given the initial segment $J$ in the fundamental region $R_F$, let $C_{[F,J]}$ be the candidate set of the leftmost tiles at the next step. Let $C^0_{F'}$ be the candidate set of the leftmost tiles at the first step of the leftmost tiling process for the fundamental region $R_{F'}$. (Recall that $C^0_{F'}$ is obtained by setting $w_X = 0$ where $X=F'$ with the notation as in Lemma \ref{LemCandMin}.) 

In order to calculate the error term $e$, we need to clarify the difference between $C^0_{F'}$ and $C_{[F,J]}$. As there are two possible cases for a tiling transition, they will be discussed in the following two lemmas, respectively. We will use the notation $\bar{C}_{[F,J]}$ and $\bar{C}^0_{F'}$ for the set $\bar{C}$ in Lemma \ref{LemCandMin} corresponding to $C_{[F,J]}$ and $C^0_{F'}$, respectively.

\begin{lemma}	\label{LemCandSetFund0}
	Suppose $[F,J]$ achieves the deducting moment with $d$ full columns and $F' = \phi([F,J])$. Then, $w \in C_{[F,J]}$ if and only if $w - d \in C^0_{F'}$.
\end{lemma}

\begin{proof}
	From the definition of the tiling transition $F' = \phi([F,J])$, it follows that $R_{[F,J]}$ and $R_{F'}$ are similar regions. Therefore, we have $w \in \bar{C}_{[F,J]}$ if and only if $w - d \in \bar{C}^0_{F'}$ by the definition of $\bar{C}_{[F,J]}$ and $\bar{C}^0_{F'}$. From the definition of the candidate set $C$ in Lemma \ref{LemCandMin}, it follows that $w \in C_{[F,J]}$ if and only if $w - d \in C^0_{F'}$. 
\end{proof}

\begin{lemma}	\label{LemCandSetFund1}
	Suppose $[F,J]$ does not achieve the deducting moment, but $F' = \pi([F,J])$ is a fundamental sequence. Then, $C_{[F,J]} \subseteq C^0_{F'}$.
\end{lemma}

\begin{proof}
	Since $F'$ is a fundamental sequence, hence by the definition of a fundamental sequence, there is an initial segment $J'$ ending with $0$ such that $\phi(J') = F'$. Let $A_{J'}$ be the essential sequence of $J'$. By Lemma \ref{LemGenSeqConnectivity}, $A_{J'}$ has connectivity. From the proof of Lemma \ref{LemValidIncreaseSeq}, we have already noted that the connectivity of $\bar{C}_{J'}$ is inherited from $A_{J'}$. ($\bar{C}_{J'}$ corresponds to $\bar{C}$ in the notation of Lemma \ref{LemCandMin} applied to $J'$.) Since $R_{F'}$ is similar to $R_{J'}$, it follows that the connectivity of $\bar{C}^0_{F'}$ is guaranteed by $\bar{C}_{J'}$, thus $C^0_{F'} = \bar{C}^0_{F'}$.
	
	 From the definition of $\bar{C}_{[F,J]}$ and $\bar{C}^0_{F'}$, it follows that $\bar{C}_{[F,J]} \subseteq \bar{C}^0_{F'}$. Note that $C_{[F,J]} \subseteq \bar{C}_{[F,J]}$ holds by the definition of $C$ in Lemma \ref{LemCandMin}. Therefore, $C_{[F,J]} \subseteq \bar{C}_{[F,J]} \subseteq \bar{C}^0_{F'} = C^0_{F'}$.
\end{proof}

The following example shows that it is possible to obtain $C_{[F,J]} \subsetneq C^0_{F'}$. For $n=4$, let $J = [5]$ and $J$ transfers $F = [1,2]$ to $F'=[1, 2, 5]$. We see that $C_{[F,J]} = \{ 3, 6, 9 \}$ and $C^0_{F'} = \{ 0, 3, 6, 9\}$. (See Figure \ref{FigTranFundamentaSeq2}.) Then, $0 \in C^0_{F'} \setminus C_{[F, J]}$. The situation $C_{[F,J]} \subsetneq C^0_{F'}$ results in the error $e$ in Equation \ref{EqNumTilingTran}.

When $[F,J]$ achieves the deducting moment, Lemma \ref{LemCandSetFund0} shows that there is a bijection between $C_{[F,J]}$ and $C^0_{F'}$. Hence, this case will not contribute to the error term $e$ in Equation \ref{EqNumTilingTran}. When $[F,J]$ does not achieve the deducting moment but $F' = \pi([F,J])$ is a fundamental sequence, Lemma \ref{LemCandSetFund1} indicates that the error term $e$ could come from an over-counting issue, and the previous example shows that the over-counting issue can really happen.

This over-counting issue can be corrected by considering the set $ C^0_{F'} \setminus C_{[F, J]} $. For every $w \in C^0_{F'} \setminus C_{[F, J]} $, let $f_w$ be the number of tiling sequences of $R_{F'}(N-|F|-|J|)$ that start with $w$. By seeking all possible tiling transitions $F' \xrightarrow{J'} F''$ where $J'$ starts from $w$, we can represent $f_w$ to be the summation of $f_{F''}(N-|F'|-|J'|)$ over all possible $F''$. We switch the sign of these terms and add them to Equation \ref{EqNumTilingTran}. 

If the over-counting issue happens again in the calculation of $f_w$, we repeat the procedure by considering the set $C^0_{F''} \setminus C_{[F', J']}$. In the whole procedure, the sign of the terms may be switched multiple times based on the inclusion-exclusion principle. In our numeric calculation, the over-counting issue happens rarely.

In the previous example, we have $C_{[F,J]} = \{ 3, 6, 9 \}$ and $C^0_{F'} = \{ 0, 3, 6, 9\}$, so $C^0_{F'} \setminus C_{[F, J]} = \{ 0 \}$. Note that $[0]$ transfers $F'=[1,2,5]$ to $F''=[2]$. Therefore, $f_{[2]}(N-4)$ is an additional negative term in the equation calculating $f_{[1,2]}(N-2)$.

We now construct the recursive system. As we will show later in Lemma \ref{LemNumFundamental}, there are $(n-1)!$ fundamental sequences. In addition, it will be shown in Lemma \ref{LemUniqueTransition} that it is sufficient to consider $n$ regions of different sizes corresponding to each fundamental sequence. Thus, define an $n!$ dimensional vector $\mathbf{f}(N) = \big( f_{F_i}(N-j) \big)^\intercal$ for $i=1,2,\cdots,(n-1)!$ and $j=0,1,\cdots,(n-1)$ such that each component $f_{F_i}(N-j)$ is the number of tilings of the fundamental region $R_{F_i}(N-j)$. We order the elements of $\mathbf{f}(N)$ first by decreasing size for fundamental regions, then for each size by the order of the sequence in $\mathbb{F}$. Note that the first component of $\mathbf{f}(N)$ is exactly the number of tilings of the $n$-by-$N$ strip. For example in the case $n=3$, the corresponding vector is 
\begin{align*}
	\mathbf{f}(N) = \Big( & f_{[0]}(N),  \ f_{[1]}(N), \ f_{[0]}(N-1), \  f_{[1]}(N-1), \  f_{[0]}(N-2), \  f_{[1]}(N-2) \Big)^\intercal	.
\end{align*}

Algorithm \ref{AlgGenFormula} in Appendix is used to obtain the transfer matrix. We report the results for $n=3$ and $n=4$. The corresponding transfer matrices $A_3$ and $A_4$ have the form declared in Theorem \ref{ThmDiffEq}. In particular, we have
$$
	A_3^{'} = 
  \begin{bmatrix}
   1 & 1 & 0 & 1  \\
   1 & 0 & 1 & 0  
  \end{bmatrix} 
$$
and 
$$
A_4^{'} =
\begin{bmatrix}
  1 & 1 & 1 & 0 & 0 & 0 & 0 & 0 & 1 & 0 & 1 & 0 & 0 & 0 & 0 & 1 & 0 & 0 \\ 
  1 & 0 & 0 & 1 & 1 & 0 & 0 & 0 & 0 & 0 & 0 & 1 & 0 & 0 & 0 & 0 & 1 & 0 \\ 
  0 & 1 & 0 & 1 & 0 & 0 & 0 & 0 & 0 & 1 & 0 & 0 & 1 & 0 & 0 & 0 & 0 & 1 \\ 
  1 & 0 & 0 & 0 & 0 & 1 & 1 & 0 & -1 & 0 & 0 & 0 & 0 & 0 & 1 & 0 & 0 & 0 \\ 
  0 & 1 & 0 & 0 & 0 & 0 & 1 & 0 & 0 & 0 & 0 & 0 & 0 & 1 & 0 & 0 & 0 & 0 \\ 
  0 & 0 & 1 & 0 & 0 & 0 & 0 & 1 & 0 & 0 & 0 & 0 & 0 & 0 & 0 & 1 & 0 & 0 \\ 
   \end{bmatrix}	.
$$

\subsection{The Generating Function}
Define the \emph{generating function} for a square matrix $A$ to be the matrix $G(x)$ with entries
$$
	G_{i,j}(x) = \sum_{N=0}^{\infty} (A^N)_{i,j} x^N	.
$$
By Theorem 4.7.2 in \citep{stanley1986enumerative}, we have
$$
	G_{i,j}(x) = \frac{(-1)^{i+j} \det(I - xA : i, j)}{\det(I-xA)}
$$
where $(B: i, j)$ denotes the matrix obtained by removing the $i$-th row and $j$-th column of $B$. 

For the number of tilings, let $g_n(x) = \sum_{N=0}^\infty f_N x^N$ be the generating function for $f_N$. Using above formula, we have 
\begin{equation} 	\label{EqGenFun}
	g_n(x) = G_{1,1}(x) = \frac{\det(I - xA_n : 1, 1)}{\det(I-xA_n)}.
\end{equation}
From Equation \ref{EqGenFun}, we are able to calculate the generating function of transfer matrix $A_n$ for small $n$. For the case $n=3$, the generating function is 
\begin{equation*}
	g_3(x) = \frac{1-x^3}{1-x-x^2-4x^3+x^6}	.
\end{equation*}
For the case $n=4$, the generating function is $g_4(x) = \frac{p(x)}{q(x)}$ where 
\begin{align*}
	p(x) = & 1 - x^2 - 13 x^4 - 2 x^5 + 5 x^6 - x^7 + 39 x^8 + 6 x^9 + 6 x^{11}  \\
 			 &	-37 x^{12} - 5 x^{14} - x^{15} + 11 x^{16} + x^{18} - x^{20}	
\end{align*}
	
and 
\begin{align*}
	q(x) = & 1 - x - 2 x^2 - 2 x^3 - 25 x^4 + 3 x^5 + 12 x^6 + 4 x^7 + 109 x^8 \\
	          & + 5 x^9 - 9 x^{10} + 7 x^{11} - 159 x^{12} - 4 x^{13} - 16 x^{14} - 7 x^{15} \\
	          & + 82 x^{16} + 10 x^{18} + x^{19} - 16 x^{20} - x^{22} + x^{24}	.
\end{align*}

\subsection{Proof of Theorem \ref{ThmDiffEq}}
\begin{figure}[H]
	\centering
	\includegraphics[scale=0.35]{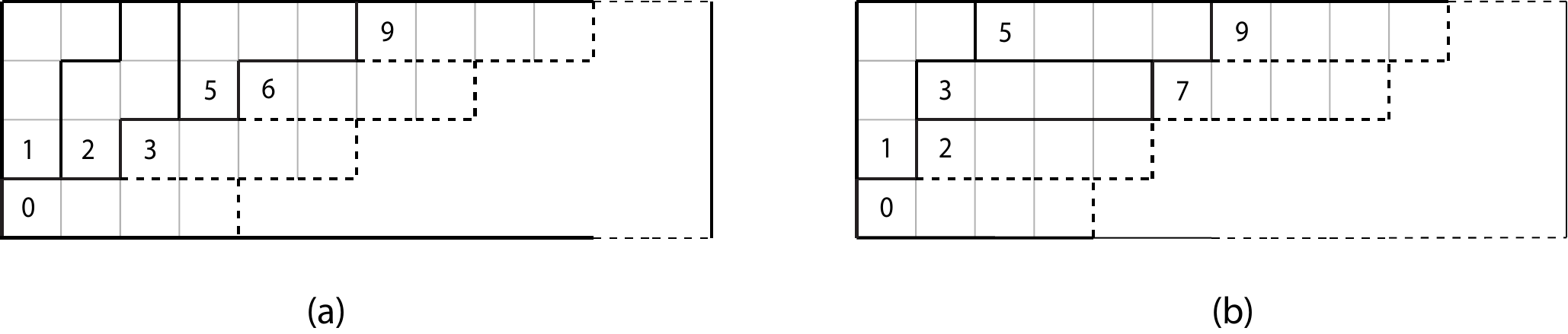}
	\caption{ $n=4$. In (a), $F = [ 1,2,5 ]$ is a fundamental sequence and $B=[ 0, 3, 6, 9 ] $ is a boundary sequence. In (b), let $J = [1, 5, 3]$ be an initial segment. The sequence $[ 0, 2, 7, 9 ]$ is not a valid boundary sequence since $7 - 2 = 5 > n$, but it is a quasi-boundary sequence. 
}
	\label{FigBoundarySequence}
\end{figure}

We start this section by describing fundamental regions using another construction. This will allow us determine the number of fundamental regions explicitly. Define a \emph{boundary sequence} $B= [ b_0, b_1, \cdots, b_{n-1} ]$ as a sequence of $n$ tiles with dimension $1 \times n$ such that $b_0 = 0$ and $b_k$ satisfies the following conditions for $k=1,\cdots,n-1$: 
\begin{itemize}
	\item[(1)] $b_{k-1}+1 \leq b_k \leq b_{k-1} + n$;
	\item[(2)] $b_k \neq b_i + ln$ for $i=0,\cdots,k-1$ and every $l \in \mathbb{N}$.
\end{itemize}
See Figure \ref{FigBoundarySequence} (a) for an example of boundary sequence $[0,3,6,9]$ and (b) for a counter example $[0,2,7,9]$. Let $\mathbb{B}$ be the set of all boundary sequences. 

Define a map $\sigma: \mathbb{F} \to \mathbb{B}$ from fundamental sequences to boundary sequences as follows. Given $F \in \mathbb{F}$, let $p$ be the left boundary of the fundamental region $R_F$, then $\sigma(F) = [ b_0, b_1, \cdots, b_{n-1} ]$ such that $b_i$ is a tile with dimension $1 \times n$ whose left boundary belongs to $p$. First, we need to show $\sigma(F) \in \mathbb{B}$. By the definition of a fundamental sequence, we have $0 \notin F$. It follows that $b_0 = 0$. Let $B' = [ b_{n-1}, \cdots, b_1, b_0 ]$. Using Lemma \ref{LemCandMin}, we can check that $[F, B']$ is a valid tiling sequence for the region $\bar{R}^F$ which is the union of $R^F$, formed by the tiles given by $F$, and $n$ horizontal tiles in rows $0,1,2, \cdots, n-1$. It follows that $\sigma(F)$ satisfies condition (1) and (2) in the definition of boundary sequence. Hence, $\sigma$ is well defined.

\begin{lemma} 	\label{LemBoundaryRigionTilable}
	Let $B$ be a boundary sequence and $R^B$ be the region of the strip on the left of the tiles of $B$. Then $R^B$ is either empty or it has  a tiling containing the tile $T_0$ with label $b_{n-1} - n$, composed of the rightmost squares at the levels $b_{n-1} - n, b_{n-1} - n +1, \cdots, b_{n-1} - 1$ in $R^B$.
\end{lemma}

\begin{pf}
	We proceed by induction on $b_{n-1}$.  The situation is trivial if $b_{n-1}$ is smallest possible, i.e. if $b_{n-1} = n-1$, since then $R^B = \emptyset$.
	
	Now let $b_{n-1} \geq n$. Then there exists $i \in \{ 1, \cdots, n-1 \}$ such that $b_{i-1} < b_{n-1} - n < b_i$. The strictness of the inequalities comes from the second part of the definition of a boundary sequence. We define
	\begin{equation*}		
		\tilde{b}_k = 
		\left\{
             \begin{aligned}
             & b_k,  & & 0 \leq k \leq i-1,  \\
             & b_{n-1} - n, & & k = i,	\\
             & b_{k-1}, & & i+1 \leq k \leq n-1.
             \end{aligned}	
		\right.
	\end{equation*}
Then $\tilde{B} = [\tilde{b}_0, \cdots, \tilde{b}_{n-1}]$	is a boundary sequence: the first condition is evident
and the second one is satisfied, because the sets of integers $\{ b_0, \cdots, b_{n-1} \}$ and $\{ \tilde{b}_0, \cdots, \tilde{b}_{n-1} \}$ agree modulo $n$. Note that $R^B$ splits into $R^{\tilde{B}}$ and $T_0$ without overlap. Using the tiling of $R^{\tilde{B}}$ given by the induction hypothesis, we obtain the required tiling of $R^B$.
\end{pf}

\begin{lemma} \label{LemBoundaryTile}
The map $\sigma$ is a bijection between the set of fundamental sequences $\mathbb{F}$ and the set of boundary sequences $\mathbb{B}$.
\end{lemma}

\begin{pf}
First, we show that $\sigma$ is an injection. For any two different fundamental sequences $F_1$ and $F_2$, their corresponding fundamental regions have different left boundaries. It follows that the corresponding boundary sequences $\sigma (F_1)$ and $\sigma (F_2)$ are different.

Now, we show that $\sigma$ is a surjection. Suppose $B = [ b_0, b_1,\cdots, b_{n-1} ]$ is a boundary sequence. By Lemma \ref{LemBoundaryRigionTilable}, let $J = [j_1, \cdots, j_k = b_{n-1} - n]$ be an initial sequence corresponding to a tiling of $R^B$. In order to prove that $\sigma$ is a surjection, it is sufficient to show that $J$ represents a fundamental sequence. Note that $\tilde{J} = J \cup [ b_{n-1}, \cdots, b_1, b_0 ]$ is a tiling sequence achieving the deducting moment with $n$ full columns. It follows that $R_{\tilde{J}}$ is a fundamental region. Note that $R_J$ and $R_{\tilde{J}}$ are similar regions since $b_i$'s are tiles with dimension $1 \times n$. By the definition of fundamental sequences, $J$ is a fundamental sequence.
\end{pf}

\begin{lemma} \label{LemMaxStep}
Running the leftmost tiling process in a strip, the deducting moment will appear in at most $L = n(n-1)/2 + 1$ steps.
\end{lemma}

\begin{pf}
Suppose $X$ is a tiling sequence at the deducting moment. Maximizing the number of elements in $X$ is equivalent to maximizing the area of the tiled region by $X$. By Lemma \ref{LemBoundaryTile}, it is equivalent to building a boundary sequence $B$ such that the area of $\bar{R}^B$ is maximized where $\bar{R}^B$ is the union of $R^B$ and the region covered by the tiles in $B$. It follows that $B$ is an arithmetic progression with initial term $0$ and common difference $n-1$. (See Figure \ref{FigBoundarySequence} (a) for an example.) Then the maximal area of $\bar{R}^B$ is $n^2(n-1)/2 + n$, which contains $L = n(n-1)/2 + 1$ tiles.
\end{pf}

Lemma \ref{LemMaxStep} gives an upper bound on the number of the steps for the appearance of the deducting moment. It guarantees that Algorithm \ref{AlgGenFundSeq} and Algorithm \ref{AlgGenFormula} are feasible. By the bijection in Lemma \ref{LemBoundaryTile}, we can also obtain all fundamental sequences using the construction of boundary sequences.

\begin{lemma} \label{LemNumFundamental}
	There are exactly $(n-1)!$ fundamental sequences.
\end{lemma}

\begin{pf}
	By Lemma \ref{LemBoundaryTile}, it is sufficient to show that there are $(n-1)!$ different boundary sequences.  Let $[ b_0, b_1, \cdots, b_{n-1} ]$ be a boundary sequence. By definition, there is only one choice for $b_0$, that is, $b_0 = 0$. Consider $b_k$ for $1 \leq k \leq n-1$. There are $n$ candidates for $b_k$ satisfying the condition (1) in the definition of boundary sequence. For these $n$ candidates of $b_k$, consider the following statement corresponding to condition (2) in the definition. For each $i$, $0 \leq i \leq k-1$, we can find an $l_i \in \mathbb{N}$ such that $b_{k-1}+1 \leq b_i + l_i n \leq b_{k-1} + n$. It follows that condition (2) will remove one of the $n$ candidates of $b_k$ for each $i$, $0 \leq i \leq k-1$. Therefore, there are $n-k$ possible values for $b_k$. Then, the number of possible choices for $[b_0, b_1, \cdots, b_{n-1}]$ is $1, n-1, n-2, \cdots, 1$. Hence, there are $(n-1)!$ possible boundary sequences. 
\end{pf}

Let $J$ be an initial segment such that $0 \notin J$, that is, the deducting moment has not been achieved. For the residual region $R_J$, let $W_J= [ w_0, w_1, \cdots, w_{n-1} ] $ be an increasing sequence of $1 \times n$ tiles such that the left boundary of these tiles coincides with the left boundary of $R_J$. We call $W_J$ the \emph{quasi-boundary sequence} of $R_J$. (For example, in Figure \ref{FigBoundarySequence} (b), $[0,2,7,9]$ is a quasi-boundary sequence of $R_{[1,5, 3]}$.) 

Since $W_J$ is adjacent to a tilable region that is covered by $J$, it follows that $W_J$ satisfies Condition (2) in the definition of boundary sequence. We observe that if $W_J$ has connectivity (satisfies Condition (1) in the definition of boundary sequence), then $W_J$ is a valid boundary sequence. From the bijection between $\mathbb{F}$ and $\mathbb{B}$ shown in Lemma \ref{LemBoundaryTile}, it follows that the  residual region $R_J$ is a fundamental region if and only if $W_J$ has connectivity. Now we are going to prove the crucial property of tiling transitions.

\begin{lemma} \label{LemUniqueTransition}
	Let $F \xrightarrow{J} F'$ be a tiling transition between fundamental sequences $F$ and $F'$ where $F$ and $F'$ are fixed. Then, $J$ is a unique decreasing sequence. Moreover, $|J| \leq n$ for every tiling transition and $|J| = n$ if and only if $J$ transits from a fundamental sequence $F$ to $F$ itself. 
\end{lemma}

\begin{figure}[H]
	\centering
	\includegraphics[scale=0.5]{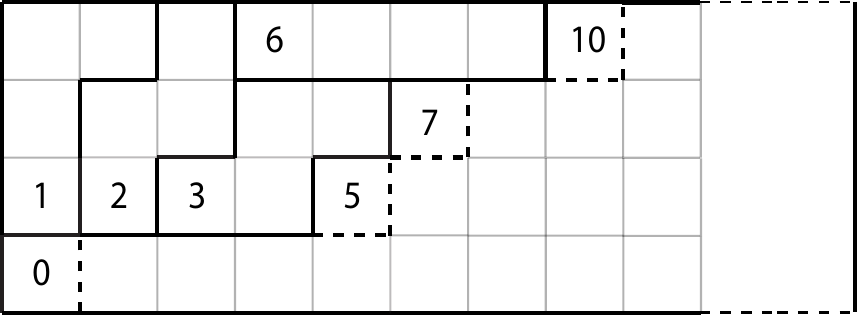}
	\caption{ $n=4$. Consider the fundamental region $R_F$ where $F = [ 1,2 ]$. Let $J^{(m)}=[ 6, 3 ] $ for $m=2$. At this step ($m=2$), the quasi-boundary sequence is $W^{(m)} = [0, 5, 7, 10]$ and the connectivity is broken between $0$ and $5$. At the previous step, $W^{(m-1)} = [0, 3,5,10]$. The formula $W^{(m)} = ( W^{(m-1)} \setminus \{ j_m \} ) \cup \{ j_m + n \}$ is satisfied. From our notation in the proof of Lemma \ref{LemUniqueTransition}, we have $W^{(m)}_* = [5, 7, 10]$, $W^{(m-1)}_* = [10]$ and $W^{(m-1)}_{**} = [5]$, which verifies the identity  $W^{(m)}_* = W^{(m-1)}_{*} \cup W^{(m-1)}_{**} \cup \{ j_m + n \}$.
}
	\label{FigTransitionSequence}
\end{figure}

\begin{pf}
By definition, $J$ is a tiling sequence from $R_F(N-l)$ to $R_{F'}(N-l')$. It is clear that $|J| = l' - l$. Let $D$ be the difference region between $R_{F'}(N-l')$ and $R_F(N-l)$ and $D$ be covered by the tiling corresponding to $J$. From \citep{sheffield2002ribbon}, it follows that the elements of $J$, which are the labels of tiles, are determined by the region $D$. For these elements of $J$, we will show that $J = [ j_1, j_2, \cdots, j_{l'-l} ]$ is a decreasing sequence. 

If $l' - l=1$, then there is only one tile $j_1$ between $F$ and $F'$. Thus, $J = [j_1]$ is uniquely determined by $F$ and $F'$. 

Suppose $l'-l > 1$. For induction, suppose $J^{(m)} = [j_1, j_2, \cdots, j_m]$ is a decreasing sequence where $m < l'-l$. Let $W^{(m)} = W_{[F, J^{(m)}]} = [w_0, \cdots, w_{n-1}]$ be the quasi-boundary sequence of the residual region $R_{[F, J^{(m)}]}$, $w_0 < w_1 < \cdots < w_{n-1}$. Since $m < l'-l$, it follows that $R_{[F, J^{(m)}]}$ is not a fundamental region by the definition of tiling transition. It follows that $W^{(m)}$ is not a valid boundary sequence, thus the connectivity of $W^{(m)}$ does not hold. We are going to prove that this implies that $j_{m+1} < j_m$.

By our induction assumption, $J^{(m)}$ is a decreasing sequence. From the definition of the leftmost tiling process, we must have a directed path $j_1 \to j_2 \to \cdots \to j_m$. It follows that $J^{(m)}$ has connectivity.

At the $m$-th step, $W^{(m)}$ can be obtained from $W^{(m-1)}$ by removing $j_m$ and including $j_m + n$, that is, $W^{(m)} = (W^{(m-1)} \setminus \{ j_m \} ) \cup \{ j_m + n \} $ where $W^{(0)}$ is the boundary sequence of $F$. It is clear that $j_m \notin W^{(m)}$ and $j_m + n \in W^{(m)}$. Let $W^{(m)}_{*} = \{ w \in W^{(m)} : w > j_m \} $ be the tail of $W^{(m)}$ where $W^{(0)}_{*} = \emptyset$,  and $W^{(m-1)}_{**} = \{ w \in W^{(m-1)}: j_m < w < j_{m-1} \}$ where $W^{(0)}_{**} = \{ w \in W^{(0)} : w > j_1 \}$. (See Figure \ref{FigTransitionSequence} for an example.) From $W^{(m)} = (W^{(m-1)} \setminus \{ j_m \} ) \cup \{ j_m + n \} $, we obtain that 
	\begin{align*}
		W^{(m)}_{*} &= \{ w \in W^{(m)} : w > j_m \}	\\
						&= \{ w \in W^{(m-1)} : w > j_m \} \cup \{ j_m + n \} 	\\
						&= \{ w \in W^{(m-1)} : w > j_{m-1} \} \cup \{ w \in W^{(m-1)} : j_m < w < j_{m-1} \} \cup \{ j_m + n \} 	\\
					     &= W^{(m-1)}_{*} \cup W^{(m-1)}_{**} \cup \{ j_m + n \}	.
	\end{align*}

We show that $W^{(m)}_{*}$ has connectivity by induction. (By default, we permute $W^{(m)}_{*}$ to be an increasing sequence.) At the first step ($m=1$), we have $W^{(1)}_{*} = W^{(0)}_{**} \cup \{ j_1 + n \}$. The sequence $W^{(0)}_{**} = \{ w \in W^{(0)} : w > j_1 \}$ has connectivity since it is the tail of the boundary sequence $W^{(0)}$. The connectivity is not broken if $j_1 + n$ is added. Hence, $W^{(1)}_*$ has connectivity. 

Suppose, for induction, that $W^{(m-1)}_*$ has connectivity. From the connectivity of $J^{(m)}$, it follows that for every $w' \in W^{(m-1)}_{**}$ we have $0 < (j_m + n) - w' < n$. Thus, the connectivity of $W^{(m-1)}_{**} \cup \{ j_m + n \}$ is proved. Note that $j_{m-1} + n \in W^{(m-1)}_*$. It follows that the connectivity between $W^{(m-1)}_{*}$ and $W^{(m-1)}_{**} \cup \{ j_m + n \}$ is guaranteed by the fact $0 < j_{m-1} - j_m < n$. Therefore, $W^{(m)}_{*} = W^{(m-1)}_{*} \cup W^{(m-1)}_{**} \cup \{ j_m + n \}$ has connectivity.

Now we consider the sequence $W^{(m)} \setminus W^{(m)}_{*}$. We see that $W^{(m)} \setminus W^{(m)}_{*}$ is the same as the head of the boundary sequence of $R_F{(N-l)}$. (For example, in Figure \ref{FigTransitionSequence}, $W^{(m)} \setminus W^{(m)}_{*} = [0]$ for $m=2$, which is unchanged as the head of the boundary sequence $W^{(0)}$.) From the definition of a boundary sequence, it follows that the connectivity of $W^{(m)} \setminus W^{(m)}_{*}$ holds. 

Let $w_i$ be the largest element of $W^{(m)}$ such that $w_i < j_m$, then $W^{(m)} \setminus W^{(m)}_{*} = [w_0, \cdots, w_i]$ and $W^{(m)}_{*} = [w_{i+1}, \cdots, w_{n-1}]$. Since the connectivity of both $W^{(m)}_{*}$ and $W^{(m)} \setminus W^{(m)}_{*}$ hold, and $W^{(m)}$ does not have connectivity, we must have $w_{i+1} - w_i > n$. (See Figure \ref{FigTransitionSequence} for an example. )

We apply Lemma \ref{LemCandMin} to $J^{(m)}$. Note that $j_m$ is the last element of $J^{(m)}$. The definition of $w_i$ ensures that $w_i$ and $w_{i+1}$ satisfy the inequalities of $\bar{C}$. It is also clear that $w_i - n, w_{i+1}-n \notin Y $ where $X = [F, J^{(m)}]$. Therefore, we have $w_i, w_{i+1} \in \bar{C}$ by definition. Since the connectivity is broken from $w_i$ to $w_{i+1}$, it follows that $w_i$ is the largest element of the candidate set $C$ of the leftmost tile at the $(m+1)$-th step of the leftmost tiling process. Therefore, $j_{m+1} \leq w_i$ since $j_{m+1} \in C$. From $w_i <  j_m$, it follows that $j_{m+1} <  j_m$. Thus, $J$ is a decreasing sequence by induction.

The maximal value of $|J|$ happens when the sequence $J$ equals $W_F$ listed in reverse order where $W_F$ is the quasi-boundary sequence of $R_F(N-l)$. (In fact, $W_F$ is a boundary sequence since $R_F(N-l)$ is a fundamental region.) In this case, we see that $R_F(N-l)$ and $R_{F'}(N-l')$ are similar, thus $F'=F$. (See Figure \ref{FigBoundarySequence} (a) for an example.) 
\end{pf}

Now we prove Theorem \ref{ThmDiffEq}.

\begin{pf}
We have constructed the recursive system for the number of tilings of fundamental regions in Section 5.1 and Section 5.2. By the uniqueness of tiling transitions proved in Lemma \ref{LemUniqueTransition}, it follows that the transfer matrix $A_n$ only contains elements $0, \pm1$. By Lemma \ref{LemNumFundamental}, there are $(n-1)!$ fundamental sequences. By Lemma \ref{LemUniqueTransition}, the size difference of the recursive system is $n$. 

Note that the difference between the candidate sets of $\tilde{C}$ based on the setup with $\tilde{w}_X = 0$ and $C$ based on the original $w_X$ may lead to an extension of tiling transitions. For every element $w \in \tilde{C} \setminus C$, we have $w < \min (C)$. Therefore, Lemma \ref{LemUniqueTransition} holds for tiling transitions including extended transitions.

The upper-right block of $A_n$ must be identity matrix $I_2$ with dimension $(n-1)! \times (n-1)!$. Therefore, the block $A_n^{'}$ has dimensions $(n-1)! \times (n-1)(n-1)!$. The identity matrix $I_1$ results from the construction of the recursive system automatically.
\end{pf}

\section{Appendix}

\subsection{Proof of Lemma \ref{LemCandMin}}
\begin{proof}
\textbf{{Proof of Necessity}.}	
We will prove $w_c \in \bar{C}$ assuming that $w_c$ is the leftmost. We first show that $w_c - n \notin Y$. Otherwise, if $w_c - n \in Y$ then the relation $w_c - n \prec w_c$ is forced. It is not possible that $w_c$ is the leftmost, contradiction.

	In order to prove the inequality $w_c \leq \max(X) + n$ in the definition of $\bar{C}$, by seeking contradiction, we suppose that $w_c -n > \max(X)$. We have already shown that $w_c - n \notin Y$. It follows that $w_c - n \in X$. However, $w_c - n > \max(X)$ contradicts the assumption that $\max(X)$ is the largest element in $X$.
	
	We prove the other inequality $w_c > w_X - n$ in the definition of $\bar{C}$ by considering the following two cases and showing that they both lead to a contradiction. 
	\begin{itemize}
		\item Case (1) $w_c = w_X - n$. In this case, the relation $w_c \prec w_X$ is forced. However, in the previous step $w_X$ was declared to be the leftmost. It contradicts the relation $w_c \prec w_X$.
		\item Case (2) $w_c < w_X - n$. In this case, $w_c$ and $w_X$ are not comparable. Consider the previous step when $w_X$ is declared to be the leftmost. One of the following two statements must be true at this step. (i) The relation $w_X \prec \cdots \prec w_c$ holds. (ii) There exists a source $v \in Y$ such that $v \prec \cdots \prec w_c$ and $w_X < v$.
		
		 For statement (i), there exists a non-empty sequence $u_1, \cdots, u_j \in Y$ such that $w_X \prec u_1 \prec \cdots \prec u_j \prec w_c$. (This sequence cannot include any of vertices in $X$ since then one of these vertices would not be a source at the moment it is chosen.) Hence, $w_c$ is not a source in $Y$ and cannot be leftmost. For statement (ii), we see $w_c$ cannot be the leftmost in $Y$ since $v \prec \cdots \prec w_c$ and $v \in Y$. Both statements (i) and (ii) contradict that $w_c$ is the leftmost.
	\end{itemize}

	To finish the proof of the necessity, we show that $w_c \in C$ assuming the existence of $i$ in the definition of $C$. (If $i$ does not exist, then $C=\bar{C}$ and there is nothing to prove.) Let  $C_0 = \{ w \in Y : w \leq w_i \}$ and $D_0 = \{ w \in Y: w_i < w \leq w_i+n \}$. We claim that $w-n \in C_0$ for every $w \in D_0$. 
	
	Now we prove the claim. For every $w \in D_0$, it is obvious that $w-n \leq w_i$, thus it is sufficient to show that $w-n \in Y$. In order to show $w-n \in Y$, we suppose $w-n \in X$ (or $w-n < 0$) seeking contradiction. From $w-n \in X$ (or $w-n < 0$), it follows that $w-n \leq \max(X)$, so $w \leq \max(X) + n$. 
	
	Since  $w_i \in C \subseteq \bar{C}$, it follows that $w_X - n < w_i$ by the definition of $\bar{C}$. (Recall that $\bar{C} = \{ w \in Y :  w_X - n < w \leq \max(X) + n \text{ and } w - n \notin Y\}$.) For $w \in D_0$, we have $w_i < w$. Hence, $w_X - n < w_i < w$. 
	
	From the inequality $w_X-n < w \leq \max(X) + n$ shown above and the assumption $w-n \in X$ (or $w-n < 0$), it follows that $w \in \bar{C}$. Since $w \in D_0$, we have $w_i < w \leq w_i + n$. However, this inequality $w_i < w \leq w_i + n$ together with $w \in \bar{C}$ contradicts the definition of the index $i$, since $w > w_i +n$ must hold for every $w \in \bar{C}$ and $w > w_i$ by the definition of $i$. Therefore, it is not possible that $w-n \in X$ (or $w-n < 0$), thus $w-n \in Y$. Hence we proved the claim.
	
	For any tiling, let $G_{C_0}$ be the corresponding graph restricted to $C_0$, and $w_s$ be a source of $G_{C_0}$. Let $E_0$ be the set of all edges $uv$ where $u \in C_0$ and $v \in D_0$. By the definition of $C_0$ and $D_0$, $E_0$ is the cut-set of the cut $(C_0, Y \setminus C_0)$ for the graph $G_Y$. 
	
	For any $w \in Y \setminus C_0$, we show that there is no directed path from $w$ to $w_s$. In order to seek contradiction, suppose that there is a directed path $p$ from $w$ to $w_s$. Since $w_s$ is a source of $G_{C_0}$ and $E_0$ is the cut-set, it follows that $p$ must contain a directed edge $\tilde{w} \to w_s$ where $\tilde{w} \in D_0$. 
	
	From the claim we have already proved, it follows that $\tilde{w} - n \in C_0$ for $\tilde{w} \in D_0$. Since $\tilde{w} \to w_s$ is a  directed edge, and $\tilde{w} > w_s$ for $\tilde{w} \in D_0$ and $w_s \in C_0$, it follows that $\tilde{w} - n$ and $w_s$ are comparable. We have $w_s \to \tilde{w} - n$, since $w_s$ is a source and $\tilde{w} - n \in C_0$. Then, $w_s \to \tilde{w} - n \to \tilde{w} \to w_s$ generates a directed cycle, contradiction. Therefore, there is no directed path from $w$ to $w_s$ where $w \in Y \setminus C_0$.
	
	For any $w \in Y \setminus C_0$, it is clear that $w > w_i$, and then $w > w_s$. We have already shown that there is no directed path from $w$ to $w_s$ for any tiling. It follows that $w$ cannot be a candidate of the leftmost tile by definition. Thus, for every candidate $w_c$ the inequality $w_c \leq w_i$ holds, then $w_c \in C$.

\textbf{{Proof of Sufficiency}.}
	Let $\tilde{C} = \{ w \in Y: w \leq w_i \text{ and } w - n \notin Y \}$ where $w_i$ is the largest element in $C$. So, $C$ is a subset of $\tilde{C}$ obtained by adding the requirement $w > w_X - n$. It is clear that $C \subseteq \tilde{C}$, and $\tilde{C} = C$ at the first step of the leftmost tiling process. We order $\tilde{C}$ as an increasing sequence $[ w_m ]_{m=-j}^i$ where $w_m \in C$ for $1 \leq m \leq i$. We claim that $\tilde{C}$ has connectivity, that is, $w_{m+1} - w_m < n$ for $m = -j, \cdots, i-1$. (We will prove the claim later.)
	
	Let $w_c \in C \subseteq \tilde{C}$ be the tile that will be declared to be the leftmost tile. We build an acyclic orientation $\alpha_{\tilde{C}}$ of $G_{\tilde{C}}$ by setting $w_s \to w_t$ for all $w_c \leq w_s < w_t$, $w_s \to w_t$ for all $w_s < w_c < w_t$, and $w_t \to w_s$ for all $w_s < w_t \leq w_c$ where $w_s$ and $w_t$ are comparable in $\tilde{C}$. Note that there are no forced edges in $G_{\tilde{C}}$, hence $\alpha_{\tilde{C}}$ is possible and acyclic. Let $\tilde{C} = [v_{j}, \cdots, v_1, w_c, u_1, \cdots, u_k]$. From our claim that $\tilde{C}$ has connectivity, we have two directed paths $w_c \to u_1 \to \cdots \to u_k$ and $w_c \to v_1 \to \cdots \to v_j$ in $\alpha_{\tilde{C}}$. It follows that $w_c$ is the unique source of $G_{\tilde{C}}$ for the acyclic orientation $\alpha_{\tilde{C}}$.
	
	Now we extend $\alpha_{\tilde{C}}$ to $G_Y$ by setting $u \to v$ for all free edges $uv$ where $u \in \tilde{C}$ and $v \in Y \setminus \tilde{C}$, and $s \to t$ if $s < t$ for all free edges $st$ where $s, t \in Y \setminus \tilde{C}$. We call this extended orientation $\alpha_Y$. It is already known that $\alpha_{\tilde{C}}$ is acyclic, and it is clear that there is no directed cycle in $G_{Y \setminus \tilde{C}}$ for $\alpha_Y$ by the construction. Moreover, there is no edge directed into $\tilde{C}$ for $\alpha_Y$ since all forced edges directed into $\tilde{C}$ have been ruled out by the fact that $w-n \notin Y$ for every $w \in \tilde{C}$. Therefore, the orientation $\alpha_Y$ is acyclic.
	
	Consider the acyclic orientation $\alpha_Y$. By construction, $w_c$ is the unique source of $G_{\tilde{C}}$. Since there is no edge directed into $\tilde{C}$, it follows that $w_c$ is a source of $G_Y$. For any other source $w_s$ of $G_Y$, if $w_s \notin \tilde{C}$, then $w_c < w_s$ by the definition of $\tilde{C}$. Thus, $w_c$ is the leftmost tile for $\alpha_Y$. In summary, for an arbitrary tile $w_c \in C$, we can build an acyclic orientation $\alpha_Y$ such that $w_c$ is the leftmost tile. This proves the sufficiency.
	
	To complete the proof, we now prove the claim that $\tilde{C}$ has connectivity. First, the connectivity holds for $C$ by the definition of the index $i$, thus the connectivity holds for $\tilde{C}$ at the first step when $\tilde{C} = C$. We will apply induction and use the notation $\cdot^{(l)}$ to represent an object $\cdot$ at the $(l)$-th step of the leftmost tiling process.
	
	For induction, suppose the connectivity holds for $\tilde{C}^{(l)}$ at the $l$-th step. Recall that $\tilde{C}^{(l)} = \{ w \in Y^{(l)}: w \leq w_i^{(l)} \text{ and } w - n \notin Y^{(l)} \}$. Let $w_c^{(l)}$ be the leftmost tile declared at the $l$-th step. At the $(l+1)$-th step, we have $w_{X^{(l+1)}} = w_c^{(l)}$, which is the final element in $X^{(l+1)}$ declared at the $l$-th step. 
	
	Let $D^{(l)} := \{ w \in Y^{(l)}: w - n \notin Y^{(l)} \}$, $A^{(l)} := \{ w \in D^{(l)}: w < w_c^{(l)} \}$ and $A'^{(l+1)} := \{ w \in D^{(l+1)}: w < w_c^{(l)} \}$. Note that $A^{(l)} \subseteq \tilde{C}^{(l)} \subseteq D^{(l)}$. By definition, $D^{(l+1)}$ differs from $D^{(l)}$ by excluding $w_c^{(l)}$ and including $w_c^{(l)} + n$, that is, $D^{(l)} \setminus D^{(l+1)} = \{ w_c^{(l)} \}$ and $D^{(l+1)} \setminus D^{(l)} = \{ w_c^{(l)} + n \}$. It follows that $A^{(l)} = A'^{(l+1)}$. In the following proof, let $A = A^{(l)} = A'^{(l+1)}$ for notation convenience and order $A$ by increasing labels. Note that $A$ is an initial segment of $\tilde{C}^{(l)}$ by the definition of $A^{(l)}$. Since $\tilde{C}^{(l)}$ has connectivity by the inductive assumption, it follows that $A$ has connectivity. 
	
	If $A = A'^{(l+1)} = \emptyset$, then by the definition of $A'^{(l+1)}$, we must have $w > w_c^{(l)}$ for every $w \in D^{(l+1)}$. It is clear that $\tilde{C}^{(l+1)} \subseteq D^{(l+1)}$, so  $w > w_c^{(l)}$ for every $w \in \tilde{C}^{(l+1)}$. It follows that every $w \in \tilde{C}^{(l+1)}$ satisfies the condition $w > w_{X^{(l+1)}} - n$ since $w_{X^{(l+1)}} = w_c^{(l)}$. From $C^{(l+1)} \subseteq \tilde{C}^{(l+1)}$ and $\tilde{C}^{(l+1)} \setminus C^{(l+1)} = \{w \in \tilde{C}^{(l+1)}: w \leq w_{X^{(l+1)}} - n \} $, we have $\tilde{C}^{(l+1)} \setminus C^{(l+1)} = \emptyset$, and therefore $\tilde{C}^{(l+1)} = C^{(l+1)}$. Therefore, the connectivity of $\tilde{C}^{(l+1)}$ holds by the connectivity of $C^{(l+1)}$. 
	
	Let $A \neq \emptyset$. At the $(l+1)$-th step, let  $w_i^{(l+1)}$ be the largest element of $C^{(l+1)}$. We consider two alternatives depending on whether $w_i^{(l+1)} \in A$ or not. If $w_i^{(l+1)} \notin A$, then $w_i^{(l+1)} \in D^{(l+1)} \setminus A'^{(l+1)}$. By the definition of $A'^{(l+1)}$, it follows that $w_i^{(l+1)} > w_c^{(l)}$, thus $w < w_i^{(l+1)}$ for every $w \in A$. 
	
	If $w_i^{(l+1)} \in A$, we show that $w \leq w_i^{(l+1)}$ for every $w \in A$. In order to seek contradiction, suppose there is $w_{i'}^{(l+1)} \in A$ such that $w_{i'}^{(l+1)} > w_i^{(l+1)}$. Since $A$ has connectivity, we can further assume that $w_{i'}^{(l+1)} < w_i^{(l+1)} + n$. We show that $w_{i'}^{(l+1)} \in \bar{C}^{(l+1)}$ by checking the inequalities in the definition of $\bar{C}^{(l+1)}$ in the following paragraph. (Recall that $\bar{C}^{(l+1)} = \{ w \in Y^{(l+1)} :  w_{X^{(l+1)}} - n < w \leq \max(X^{(l+1)}) + n \text{ and } w - n \notin Y^{(l+1)} \}$.) 
	
	First, since $w_i^{(l+1)} \in C^{(l+1)}$, it follows that $w_i^{(l+1)} > w_{X^{(l+1)}} - n$ by the definition of $C^{(l+1)}$, then $w_{i'}^{(l+1)} > w_i^{(l+1)} > w_{X^{(l+1)}} - n$. Secondly, since $w_c^{(l)} \in C^{(l)}$, it follows that $w_c^{(l)} \leq \max(X^{(l)}) + n$. Since $w_{i'}^{(l+1)} \in A$, then by the definition of $A$ we have $w_{i'}^{(l+1)} < w_c^{(l)}$. Thus, we have $w_{i'}^{(l+1)} < w_c^{(l)} \leq \max(X^{(l)}) + n \leq \max(X^{(l+1)}) + n $. Therefore, $w_{i'}^{(l+1)}$ satisfies the inequalities in the definition of $\bar{C}^{(l+1)}$.
	
	From $w_{i'}^{(l+1)} \in \bar{C}^{(l+1)}$ and the inequality $w_i^{(l+1)} < w_{i'}^{(l+1)} < w_i^{(l+1)} + n$, we have a contradiction with the assumption that $w_i^{(l+1)}$ is the largest element of $C^{(l+1)}$. Hence, it is not possible that there is $w_{i'}^{(l+1)} \in A$ such that $w_{i'}^{(l+1)} > w_i^{(l+1)}$. Thus, $w \leq w_i^{(l+1)}$ for every $w \in A$.  
		
	For both alternatives of $w_i^{(l+1)} \in A$ and $w_i^{(l+1)} \notin A$, we have shown that $w \leq w_i^{(l+1)}$ for every $w \in A$. It follows that $A \subseteq \tilde{C}^{(l+1)}$ by the definition of $\tilde{C}^{(l+1)}$. Let $A$ and $\tilde{C}^{(l+1)}$ be ordered by increasing labels, then $A$ is an initial segment of $\tilde{C}^{(l+1)}$. 
	
	We consider the following two cases by comparing the elements of $A$ with $w_c^{(l)} - n$.
	
	\begin{itemize}
		\item Case (1) $w > w_c^{(l)} - n$ for every $ w \in A$. Since $A$ is an initial segment of $\tilde{C}^{(l+1)}$, it follows that $w > w_c^{(l)} - n$ for every $ w \in \tilde{C}^{(l+1)}$. By the definition of $\tilde{C}^{(l+1)}$ and $C^{(l+1)}$, we have $\tilde{C}^{(l+1)} \setminus C^{(l+1)} = \{ w \in \tilde{C}^{(l+1)}: w \leq w_{X^{(l+1)}} - n \}$. Since $w_{X^{(l+1)}} = w_c^{(l)}$, it follows that $\tilde{C}^{(l+1)} \setminus C^{(l+1)} = \emptyset$, thus $\tilde{C}^{(l+1)} = C^{(l+1)}$. Hence, the connectivity of $\tilde{C}^{(l+1)}$ holds by the connectivity of $C^{(l+1)}$.

		\item Case (2) There is at least one element $w_{m_0}^{(l)} \in A$ such that $w_{m_0}^{(l)} < w_c^{(l)} - n$. We have $w_{m_0}^{(l)}, w_c^{(l)} \in \tilde{C}^{(l)}$ by the definition of $\tilde{C}^{(l)}$. From the connectivity of $\tilde{C}^{(l)}$ ensured by the inductive assumption, it follows that there is an element $w_{m_1}^{(l)} \in \tilde{C}^{(l)}$ such that $w_{m_0}^{(l)} < w_c^{(l)} - n < w_{m_1}^{(l)} < w_c^{(l)}$. From $w_{m_1}^{(l)} \in \tilde{C}^{(l)}$ and the inequality $w_{m_1}^{(l)} < w_c^{(l)}$, we have $w_{m_1}^{(l)} \in A$ by the definition of $A^{(l)}$. It follows that $w_{m_1}^{(l)} \in \tilde{C}^{(l+1)}$ since $A = A'^{(l+1)} \subseteq \tilde{C}^{(l+1)}$.  
		
		Let $\tilde{C}_1^{(l+1)} = \{ w \in \tilde{C}^{(l+1)}: w \leq w_{m_1}^{(l)} \}$ and $\tilde{C}_2^{(l+1)} = \{ w \in \tilde{C}^{(l+1)}: w \geq w_{m_1}^{(l)} \}$. We have $\tilde{C}_1^{(l+1)}$ is an initial segment of $A$ since $w_{m_1}^{(l)} < w_c^{(l)}$, thus the connectivity of $\tilde{C}_1^{(l+1)}$ is preserved from $A$. We also have $\tilde{C}_2^{(l+1)}$ is a tail of $C^{(l+1)}$ since $w_{m_1}^{(l)} > w_c^{(l)} - n = w_{X^{(l+1)}} - n$, thus the connectivity of $\tilde{C}_2^{(l+1)}$ is preserved from $C^{(l+1)}$. It is clear that $\tilde{C}^{(l+1)} = \tilde{C}_1^{(l+1)} \cup \tilde{C}_2^{(l+1)}$ and $\tilde{C}_1^{(l+1)} \cap \tilde{C}_2^{(l+1)} = \{ w_{m_1} \}$. Therefore, the connectivity of $\tilde{C}^{(l+1)}$ holds.	
	\end{itemize}
\end{proof}

\subsection{Algorithms}

\begin{algorithm}[H] 	\label{AlgGenFundSeq}
\SetAlgoLined
\KwResult{Output the set $\mathbb{F}$ of all fundamental sequences. }

	Initialize $\mathbb{J} = \{ [ i ] : 0 \leq i \leq n-1 \}$, and $\mathbb{F}$ to be empty. \\
	\While{$\mathbb{J}$ is not empty}{
		\For{each $J \in \mathbb{J}$}{
			Generate the candidate set $C_J$ for the leftmost tile according to Lemma \ref{LemCandMin}. \\
			\For{each $w \in C_J$}{
				Let $J_w = [J, w]$ by putting $w$ consecutive to $J$.		\\
				Append $J_w$ into $\mathbb{J}$.
			}
			Delete $J$ from $\mathbb{J}$.
 		 }
		\For{each $J \in \mathbb{J}$}{
			\If{$0 \in J$}{
				Simplify $J$ to be a fundamental sequence $F$. That is, $F = \phi(J)$.	\\
				Append $F$ into $\mathbb{F}$ without duplication. 	\\
				Delete $J$ from $\mathbb{J}$.
			}
		}
	 }
	 Output $\mathbb{F}$.
\caption{Generate Fundamental Sequences}
\end{algorithm}

\vspace{10pt}
\textbf{Remark.} For every $J \in \mathbb{J}$, it is sufficient to store $J$ as an increasing sequence $\pi(J)$ and the last element of $J$. This information will be enough to generate candidate set of the leftmost tiles by Lemma \ref{LemCandMin}. Then, the duplication in $\mathbb{J}$ can be removed. For example, we can store both $[1,2,5]$ and $[2,1,5]$ as a two-tuple $( [1,2,5], 5)$.

\newpage
\begin{algorithm}[H] 	\label{AlgGenFormula}
\SetAlgoLined
\KwResult{Output the transfer matrix $A$. }

	Initialize $A$ to be an $n!$-by-$n!$ matrix with all $0$ elements. \\
	\For{each fundamental region $R_F(N)$ where $F \in \mathbb{F}$}{
		Generate the candidate set $C_F$ for the leftmost tiles of $R_F{(N)}$. \\
		Initialize $\mathbb{J} = \{ [w]: w \in C_F \}$ to be the set of all sequences. 	\\
		Set $signPower = 0$ for each sequence in $\mathbb{J}$.	\\
		\While{$\mathbb{J}$ is not empty}{
			Let $\mathbb{J}'$ be an empty set.	\\
			\For{each $J \in \mathbb{J}$}{
				\If{$0 \in J$}{
					Let $F' = \phi([F,J])$. (Then $F' \in \mathbb{F}$.)		\\
					Find the index $(i, j)$ of $A$ corresponding to $f_F(N)$ and $f_{F'}(N-|J|)$.	\\
					Assign $A_{i,j} = (-1)^{signPower}$. 	\\
					Delete $J$ from $\mathbb{J}$.
				}
				\If{$\pi([F,J]) = F' \in \mathbb{F}$}{ 
					Find the index $(i, j)$ of $A$ corresponding to $f_F(N)$ and $f_{F'}(N - |J|)$.	\\
					Assign $A_{i,j} = (-1)^{signPower}$. 	\\
					Generate candidate set $C_{[F,J]}$ for $J$ by Lemma \ref{LemCandMin}. 	\\
					Generate candidate set $C_{F'}$ for $R_{F'}(N - |J|)$ by Lemma \ref{LemCandMin}. 	\\
					Let $C_{diff} = C_{F'} \setminus C_{[F,J]}$.		\\
					\If{$C_{diff}$ is not empty}{
						\For{each $w \in C_{diff}$}{
							Let $J' = [J, w]$ and $signPower = signPower + 1$ for $J'$.	\\
							Add $J'$ to the set $\mathbb{J}'$.
						}
					}
					Delete $J$ from $\mathbb{J}$.
				}
			}
			\If{$\mathbb{J}$ is not empty}{
				Update all sequences in $\mathbb{J}$ by Lemma \ref{LemCandMin} with $signPower$ inherited.  \\
			}
			Let $\mathbb{J} = \mathbb{J} \cup \mathbb{J}'$.
 		 }
	 }
	 Output $A$.
\caption{Generate transfer Matrix}
\end{algorithm}

\bibliography{Bibliography.bib}
\bibliographystyle{asa.bst}

\end{document}